# STABLE STATIONARY PROCESSES RELATED TO CYCLIC FLOWS[1]

By Vladas Pipiras and Murad S. Taqqu

*Boston University*

We study stationary stable processes related to periodic and cyclic flows in the sense of Rosiński [*Ann. Probab.* **23** (1995) 1163–1187]. These processes are not ergodic. We provide their canonical representations, consider examples and show how to identify them among general stationary stable processes. We conclude with the unique decomposition in distribution of stationary stable processes into the sum of four major independent components: 1. A mixed moving average component. 2. A harmonizable (or "trivial") component. 3. A cyclic component 4. A component which is different from these.

**1. Introduction.** Consider a symmetric $\alpha$-stable (S$\alpha$S, for in short), $\alpha \in (0,2)$, *stationary* process $\{X_\alpha(t)\}_{t \in T}$ that has an integral representation

$$\{X_\alpha(t)\}_{t \in T} \stackrel{d}{=} \left\{ \int_S f_t(s) M_\alpha(ds) \right\}_{t \in T}, \tag{1.1}$$

where $\stackrel{d}{=}$ stands for equality in the sense of the finite-dimensional distributions. Here, $T = \mathbb{Z}$ or $T = \mathbb{R}$, $(S, \mathcal{S}, \nu)$ is a standard Lebesgue space (see Appendix A for a precise definition),

$$\{f_t\}_{t \in T} \subset L^\alpha(S, \mathcal{S}, \nu)$$

is a collection of deterministic functions such that the map $f_t(s): T \times S \mapsto \mathbb{R}$ or $\mathbb{C}$ is measurable and $M_\alpha$ is, respectively, either a real-valued or a complex-valued rotationally invariant S$\alpha$S random measure on $(S, \mathcal{S})$ with the control measure $\nu$. (Rotationally invariant means that $e^{i\gamma} M_\alpha \stackrel{d}{=} M_\alpha$ for any real

Received April 2002; revised September 2003.
[1]Supported in part by NSF Grants DMS-01-02410 and ANI-98-05623 at Boston University.

angle $\gamma$.) The process $X_\alpha$ is real-valued if the random measure $M_\alpha$ and the functions $f_t$ are real-valued; it is complex-valued if the measure and the functions are complex-valued.

Relationship (1.1) then means that the characteristic function of the process $X_\alpha$ can be expressed as

$$(1.2) \quad E\exp\left\{i\sum_{k=1}^{n}\theta_k X_\alpha(t_k)\right\} = \exp\left\{-\int_S \left|\sum_{k=1}^{n}\theta_k f_{t_k}(s)\right|^\alpha \nu(ds)\right\},$$

where $\theta_k \in \mathbb{R}, t_k \in T$, in the real-valued case, and as

$$\begin{aligned}(1.3) \quad & E\exp\left\{i\sum_{k=1}^{n}\Re(\overline{\theta}_k X_\alpha(t_k))\right\} \\ &= E\exp\left\{i\sum_{k=1}^{n}(\Re(\theta_k)\Re(X_\alpha(t_k)) + \Im(\theta_k)\Im(X_\alpha(t_k)))\right\} \\ &= \exp\left\{-c_0\int_S \left|\sum_{k=1}^{n}\overline{\theta}_k f_{t_k}(s)\right|^\alpha \nu(ds)\right\},\end{aligned}$$

where $\theta_k \in \mathbb{C}, t_k \in T$, $c_0 = (2\pi)^{-1}\int_0^{2\pi} |\cos\phi|^\alpha d\phi$ and $\overline{z}$ denotes the complex conjugate of $z \in \mathbb{C}$, in the complex-valued case [see, e.g, Samorodnitsky and Taqqu (1994)]. It is known, for example, that every measurable real-valued S$\alpha$S process $X_\alpha$ has an integral representation (1.1) with, for example, $S = (0,1)$, $\mathcal{S} = \mathcal{B}(0,1)$ and $\nu =$ Lebesgue measure [see Samorodnitsky and Taqqu (1994), Theorems 13.2.1 and 9.4.2]. Finally, recall that $\{X_\alpha(t)\}_{t\in T}$ is stationary if, for all $h \in T$, the finite-dimensional distributions of the process $X_\alpha(t+h)$, $t \in T$, are identical to those of the process $X_\alpha(t)$, $t \in T$.

In a fundamental paper, Rosiński (1995) showed that a S$\alpha$S stationary process $X_\alpha$ can be related to a flow and a corresponding cocycle as in Definition 1.1. A flow is a collection of deterministic maps $\{\phi_t\}_{t\in T}$ that satisfy

$$\phi_{t_1+t_2} = \phi_{t_1} \circ \phi_{t_2}, \qquad t_1, t_2 \in T.$$

A cocycle $\{a_t\}_{t\in T}$ for the flow $\{\phi_t\}_{t\in T}$ satisfies relationship

$$a_{t_1+t_2} = a_{t_1} a_{t_2} \circ \phi_{t_1}, \qquad t_1, t_2 \in T.$$

See Appendix A for precise definitions. By support of $\{f_t\}_{t\in T}$, we mean a minimal (a.e.) set $A \in \mathcal{S}$ such that $\nu\{f_t(s) \neq 0, s \notin A\} = 0$ for every $t \in T$. The support is denoted $\mathrm{supp}\{f_t, t \in T\}$.

DEFINITION 1.1 [Rosiński (1995)]. A S$\alpha$S stationary process $X_\alpha$ that has a representation (1.1) is said to be generated by a nonsingular measurable flow $\{\phi_t\}_{t\in T}$ on $(S, \mathcal{S}, \nu)$ if, for all $t \in T$,

$$(1.4) \qquad f_t(s) = a_t(s)\left\{\frac{d(\nu \circ \phi_t)}{d\nu}(s)\right\}^{1/\alpha} f_0(\phi_t(s)) \qquad \text{a.e. } \nu(ds),$$

where $f_0 \in L^\alpha(S, \mathcal{S}, \nu)$ and $\{a_t\}_{t \in T}$ is a cocycle for the flow $\{\phi_t\}_{t \in T}$ taking values in $\{-1, 1\}$ in the real-valued case and in the unit circle $\{w: |w| = 1\}$ in the complex-valued case, and

$$(1.5) \qquad \text{supp}\{f_t, t \in T\} = S, \qquad \nu\text{-a.e.}$$

Observe that this definition is consistent with stationarity because it implies, by using the definitions of a flow and a cocycle, that

$$\int_S \left| \sum_{k=1}^n \theta_k f_{t_k+h}(s) \right|^\alpha \nu(ds)$$

$$= \int_S \left| \sum_{k=1}^n \theta_k a_{t_k+h}(s) \left\{ \frac{d(\nu \circ \phi_{t_k+h})}{d\nu}(s) \right\}^{1/\alpha} f_0(\phi_{t_k+h}(s)) \right|^\alpha \nu(ds)$$

$$= \int_S \left| \sum_{k=1}^n \theta_k a_{t_k}(\phi_h(s)) \left\{ \frac{d(\nu \circ \phi_{t_k})}{d\nu}(\phi_h(s)) \right\}^{1/\alpha} f_0(\phi_{t_k}(\phi_h(s))) \right|^\alpha (\nu \circ \phi_h)(ds)$$

$$= \int_S \left| \sum_{k=1}^n \theta_k f_{t_k}(s) \right|^\alpha \nu(ds),$$

where the last equality follows by a change of variables $\phi_h(s) \to s$ and (1.4).

Definition 1.1 relates $f_t$ to $f_0 \circ \phi_t$. By using this connection between kernels and flows, Rosiński (1995) obtained a unique decomposition in distribution of S$\alpha$S stationary processes into two independent processes

$$(1.6) \qquad X_\alpha \stackrel{d}{=} X_\alpha^D + X_\alpha^C,$$

where the process $X_\alpha^D$ is generated by a dissipative flow and the process $X_\alpha^C$ is generated by a conservative flow (see Appendix A for definitions of dissipative and conservative flows). Moreover, Rosiński showed that dissipative processes $X_\alpha^D$ have a canonical representation

$$(1.7) \qquad \int_X \int_T k(x, t+u) M_\alpha(dx, du),$$

where $(X, \mathcal{X}, \mu)$ is a standard Lebesgue space and $M_\alpha$ has the control measure $\mu(dx)\lambda(du)$ with ($\delta$ denotes a counting measure)

$$(1.8) \qquad \lambda(du) = \begin{cases} \delta_{\mathbb{Z}}(du), & \text{if } T = \mathbb{Z}, \\ du, & \text{if } T = \mathbb{R}, \end{cases}$$

and showed that conservative processes $X_\alpha^C$ can be uniquely decomposed further into two independent processes $X_\alpha^F$ and $X_\alpha^{C \setminus F}$.

In the complex-valued case, the process $X_\alpha^F$ is the *harmonizable* process

$$(1.9) \qquad X_\alpha^F(t) \stackrel{d}{=} \int_{\widehat{T}} e^{itx} N_\alpha(dx),$$

where $N_\alpha$ has a finite control measure $\eta$ on

$$(1.10) \qquad \widehat{T} = \begin{cases} \mathbb{R}, & \text{if } T = \mathbb{R}, \\ [0, 2\pi), & \text{if } T = \mathbb{Z}. \end{cases}$$

In the real-valued case, $X_\alpha^F$ is the *trivial* stationary process

$$(1.11)\ X_\alpha^F(t) \stackrel{d}{=} X_1 + n(t)X_2 = \int_{\{1\}} N_\alpha(dx) + n(t)\int_{\{2\}} N_\alpha(dx) \qquad \text{a.s. } \forall\, t \in T,$$

where $X_1$ and $X_2$ are independent S$\alpha$S random variables, $N_\alpha$ has a finite control measure $\eta$ on $\{1, 2\}$ and

$$(1.12) \qquad n(t) = \begin{cases} 0, & \text{if } T = \mathbb{R}, \\ (-1)^t, & \text{if } T = \mathbb{Z}, \end{cases}$$

that is, $X_\alpha^F(t) \stackrel{d}{=} X_1$ if $T = \mathbb{R}$ and $X_\alpha^F(t) \stackrel{d}{=} X_1 + (-1)^t X_2$ if $T = \mathbb{Z}$.

The process $X_\alpha^F$, whether harmonizable or trivial, is a stationary S$\alpha$S process generated essentially by the simplest type of conservative flows, namely, the *identity flows* defined by $\phi_t(s) = s$ for all $t \in T, s \in S$ with $S = \widehat{T}$ in the complex-valued case and $S = \{1, 2\}$ in the real case. The other process in the decomposition of $X_\alpha^C$, that is, $X_\alpha^{C\setminus F}$, is a S$\alpha$S stationary process generated by a conservative flow which does not have a harmonizable (or trivial) component, meaning that it cannot be decomposed into two independent processes one of which is either a harmonizable process in the complex-valued case or a trivial stationary process in the real-valued case. This led Rosiński (1995) to the unique decomposition in distribution of S$\alpha$S stationary processes into three independent processes

$$(1.13) \qquad X_\alpha \stackrel{d}{=} X_\alpha^D + X_\alpha^F + X_\alpha^{C\setminus F}.$$

A nice review can be found in Rosiński (1998).

In this work, we focus on S$\alpha$S stationary processes generated by *periodic* flows in the sense of Definition 1.1 and more specifically by *cyclic* flows. Periodic flows are examples of conservative flows such that any point in the space comes back to its initial position in a finite period of time. Identity flows are periodic flows with period zero. Cyclic flows are periodic flows with positive period. We will show that S$\alpha$S stationary processes generated by periodic flows have a canonical representation which is given by the sum of two terms. The first term is the harmonizable or trivial process

$$(1.14) \quad \begin{aligned} &\int_{\widehat{T}} e^{itx} N_\alpha(dx) &&\text{(complex-valued case)}, \\ &\int_{\{1\}} N_\alpha(dx) + n(t)\int_{\{2\}} N_\alpha(dx) &&\text{(real-valued case)}; \end{aligned}$$

the second term is

$$\text{(1.15)} \qquad \int_Z \int_{[0,q(z))} b_1(z)^{[v+t]_{q(z)}} g(z, \{v+t\}_{q(z)}) M_\alpha(dz, dv),$$

where, for $a > 0$ and $x \in \mathbb{R}$,

$$\text{(1.16)} \qquad [x]_a = \max\{n \in \mathbb{Z} : na \leq x\}, \qquad \{x\}_a = x - a[x]_a \geq 0,$$

$(Z, \mathcal{Z}, \sigma)$ is a standard Lebesgue space, $q(z) \in T_+$ with

$$\text{(1.17)} \qquad T_+ = \begin{cases} (0, \infty), & \text{if } T = \mathbb{R}, \\ \{2, 3, \dots\}, & \text{if } T = \mathbb{Z}, \end{cases}$$

and

$$b_1(z) \in \begin{cases} \{w : |w| = 1\}, & \text{(complex-valued case)}, \\ \{-1, 1\}, & \text{(real-valued case)}, \end{cases}$$

$$g \in L^\alpha(Z \times [0, q(\cdot)), \sigma(dz)\lambda(dv)),$$

with $\lambda(dv)$ defined in (1.8). Moreover, $M_\alpha$ and $N_\alpha$ above are *independent* S$\alpha$S random measures with the control measures $\sigma(dz)\lambda(dv)$ and $\eta(dx)$, respectively, so that (1.14) and (1.15) are independent processes. The processes represented by the sum of (1.14) and (1.15) are called *stationary periodic processes*. Observe that the term "periodic" refers to the flow and not to the sample path behavior of the process. A stationary periodic process is called a *stationary cyclic process* if it does not have a harmonizable (or trivial) component, that is, if it cannot be represented as a sum of two independent stationary processes one of which is a nondegenerate harmonizable (or trivial) process. Note that stationary cyclic processes cannot be defined by (1.15) because, for example, harmonizable or trivial processes (1.14) can also be represented by (1.15) (see Lemma 3.1).

Stationary periodic processes (1.14) and (1.15) are always generated by periodic flows because the process (1.14) is generated by an identity flow and the process (1.15) is generated by a cyclic flow (see Theorem 3.1). We show in Theorem 3.2 that if representation (1.14) and (1.15) are *minimal*, that is, if there is no redundancy in the representation (minimal representations are defined in Appendix B), then a stationary periodic (cyclic, resp.) process can only be generated by a periodic (cyclic, resp.) flow.

If the representation is not minimal, stationary periodic processes (1.14) and (1.15) may also be generated by flows that are not periodic (see Example 3.1) and stationary cyclic processes may also be generated by flows that are not cyclic. To determine, therefore, whether a given stationary stable processes is a stationary periodic or cyclic process, it is in general not enough to examine whether the underlying flow is periodic or cyclic. There is, however, an alternate criterion that can be used to identify stationary periodic and cyclic processes. This criterion is based on the structure of the kernel function $f_t$ in (1.1) (see Theorems 4.1 and 5.1). Thus while flows have a physical interpretation, the identification criterion, which is based on the kernel, has the advantage that it can be used whether the representation is minimal or not. An analogous approach was followed by Rosiński (1995) in the case of harmonizable (or trivial) processes (1.9) [or (1.11)], typically associated with identity flows.

Our goal then is to identify stationary periodic (cyclic, resp.) processes among general S$\alpha$S stationary processes, namely, to be able to conclude that a given S$\alpha$S stationary process is a stationary periodic (cyclic, resp.) process, either by using flows in the case of a minimal representation or by applying the identification criterion mentioned above.

The identification criterion provides also a decomposition of S$\alpha$S stationary processes which is more refined than the decomposition (1.13) obtained by Rosiński (1995). More precisely, we will show (see Theorem 6.1) that the "third kind" process $X_\alpha^{C \setminus F}$ in (1.13) can be further uniquely decomposed into two independent processes

$$(1.18) \qquad X_\alpha^{C \setminus F} \stackrel{d}{=} X_\alpha^L + X_\alpha^{C \setminus P},$$

where $X_\alpha^L$ is a stationary cyclic process and $X_\alpha^{C \setminus P}$ is a S$\alpha$S stationary process generated by a conservative flow, that is without a periodic component.

A simple example of a real-valued S$\alpha$S stationary cyclic process with $T = \mathbb{R}$ is the real part of a harmonizable process (1.9),

$$(1.19) \quad \begin{aligned} \Re \int_{\mathbb{R}} e^{itx} M_\alpha(dx) &\stackrel{d}{=} c \int_{\mathbb{R}} \int_0^{2\pi} \cos(v + zt) M_\alpha(dz, dv) \\ &\stackrel{d}{=} c \int_{\mathbb{R}} \int_0^{2\pi/|z|} \cos\left(z\{w + t\}_{2\pi/|z|}\right) M_\alpha(dz, dw), \end{aligned}$$

that is, the process (1.15) with $b_1(z) = 1$, $q(z) = 2\pi/|z|$ and $g(z, u) = \cos(zu)$ (see Example 3.2). We show in Example 5.1 that the process (1.19) is indeed a stationary cyclic process, that is, an example of a process $X_\alpha^L$ in the decomposition (1.18). An example of the process $X_\alpha^{C \setminus P}$ is the stationary sub-Gaussian process (see Example 6.1).

This paper is organized as follows. In Section 2, we prove some results on periodic and cyclic flows that are used in the sequel. In Section 3, we show that stationary S$\alpha$S processes generated by periodic flows have a canonical representation given by the sum of (1.14) and (1.15). In Section 4, we provide a criterion to identify stationary periodic processes among general S$\alpha$S stationary processes. In Section 5, we do this for stationary cyclic processes. A further decomposition of stationary S$\alpha$S processes is established in Section 6. Finally, in Appendix A, we collect some basic facts related to flows and, in Appendix B, we recall the definition of minimal integral representations for stable processes.

**2. Periodic and cyclic flows.** Let $\{\phi_t\}_{t \in T}$ be a measurable flow on a standard Lebesgue space $(S, \mathcal{S}, \nu)$, where $T = \mathbb{Z}$ or $T = \mathbb{R}$ (see Appendix A). Let also

(2.1) $$P := \{s : \exists p = p(s) \in T \setminus \{0\} : \phi_p(s) = s\},$$

(2.2) $$F := \{s : \phi_t(s) = s \text{ for all } t \in T\},$$

(2.3) $$L := P \setminus F$$

be the *periodic*, *fixed* and *cyclic* points of the flow $\{\phi_t\}_{t \in T}$, respectively.

DEFINITION 2.1. A measurable flow $\{\phi_t\}_{t \in T}$ on $(S, \mathcal{S}, \nu)$ is *periodic* if $S = P$ $\nu$-a.e., is *identity* if $S = F$ $\nu$-a.e., and is *cyclic* if $S = L$ $\nu$-a.e.

We say henceforth that a set $A \subset S$ or a map $f$ on $S$ is $\nu$-measurable if it is measurable with respect to a measure $\nu$, that is, measurable with respect to the completion of the Borel sets under that measure.

LEMMA 2.1. *The set $F$ in* (2.2) *is (Borel) measurable and the sets $P$ in* (2.1) *and $L$ in* (2.3) *are $\nu$-measurable.*

PROOF. Since the proof of the lemma is elementary when $T = \mathbb{Z}$, we consider only the case $T = \mathbb{R}$. To show that the set $F$ is measurable, we first show that $F' := \{s : \phi_t(s) = s \text{ a.e. } dt\}$ satisfies $F' = F$. Indeed, if $s \in F'$, then, by definition, $\tau := \{t : \phi_t(s) = s\} = \mathbb{R}$ a.e. Observe that $\tau$ is an additive group of $\mathbb{R}$ [if $t_1, t_2 \in \tau$, then $t_1 + t_2 \in \tau$, because $\phi_{t_1+t_2}(s) = \phi_{t_1}(\phi_{t_2}(s)) = \phi_{t_1}(s) = s$] and hence by Corollary 1.1.4 in Bingham, Goldie and Teugels (1987), we have $\tau = \mathbb{R}$ and hence $F' = F$. Then, $F = \{s : h(s) = 0\}$, where $h(s) = \int_\mathbb{R} \mathbb{1}_{\{\phi_t(s) \neq s\}}(t, s) \, dt$. Since the function $h(\cdot)$ is measurable by the Fubini theorem, the set $F$ is measurable as well [use Theorem A in Halmos (1950), page 143]. To prove that the set $P = \{s : \exists p = p(s) \neq 0 : \phi_p(s) = s\}$ is $\nu$-measurable, consider the measurable set $\widetilde{P} = \{(s, p) : \phi_p(s) = s, p \neq 0\}$. Observe that $P = \text{proj}_S\{\widetilde{P}\} := \{s : \exists p : (s, p) \in \widetilde{P}\}$. The $\nu$-measurability of $P$ follows from Lemma 4.2. The set $L$ is $\nu$-measurable because $L = P \setminus F$. □

We use in the sequel the following alternative definition of a cyclic flow, which is equivalent to Definition 2.1 by Theorem 2.1.

DEFINITION 2.2. A measurable flow $\{\phi_t\}_{t \in T}$ on $(S, \mathcal{S}, \nu)$ is *cyclic* if it is null isomorphic (mod 0) to the flow

(2.4) $$\widetilde{\phi}_t(z, v) = (z, \{v + t\}_{q(z)})$$

on $(Z \times [0, q(\cdot)), \mathcal{Z} \times \mathcal{B}([0, q(\cdot))), \sigma(dz)\lambda(dv))$, where $q(z) \in T_+$ a.e. is some measurable function [see also the notation (1.8), (1.16) and (1.17)].

The $\sigma$-field $\mathcal{Z} \times \mathcal{B}([0, q(\cdot)))$ in Definition 2.2 is defined as the restriction of the $\sigma$-field $\mathcal{Z} \times \mathcal{B}(\mathbb{R})$ to the set $Z \times [0, q(\cdot))$. Null isomorphic (mod 0) in Definition 2.2 means that there are two null sets $N \subset S$ and $\widetilde{N} \subset Z \times [0, q(\cdot))$, and a Borel measurable, one-to-one, onto and nonsingular map with a measurable inverse (a so-called null isomorphism) $\Phi : Z \times [0, q(\cdot)) \setminus \widetilde{N} \mapsto S \setminus N$ such that

$$(2.5) \qquad \phi_t(\Phi(z,v)) = \Phi(\widetilde{\phi}_t(z,v))$$

for all $t \in T$ and $(z,v) \in Z \times [0, q(\cdot)) \setminus \widetilde{N}$. The null sets $N$ and $\widetilde{N}$ are required to be invariant under the flows $\phi_t$ and $\widetilde{\phi}_t$, respectively. We can view (2.4) in two ways. At each $z$ on the horizontal axis, we climb to *height* $q(z)$ before falling back on the horizontal axis to the same starting points. Equivalently, we can view the space $(z, v)$ shaped as a torus where at each point $z$ on the "grand circle," we rotate around the "small circle" of length $q(z)$.

EXAMPLE 2.1. A collection of maps

$$\phi_t(s) = e^{it\theta} s, \qquad t \in \mathbb{R},$$

with some fixed $\theta > 0$, is a measurable flow on the unit circle $\{s : |s| = 1\}$. The flow $\{\phi_t\}_{t \in \mathbb{R}}$ is cyclic since each point of the space is not fixed and comes back to its initial position in a finite time (Definition 2.1) or since it is isomorphic to the flow $\widetilde{\phi}_t(v) = \{v + t\}_{2\pi\theta^{-1}}$, $t \in \mathbb{R}$, $v \in [0, 2\pi\theta^{-1})$ (Definition 2.2). This corresponds to representing motion on the circle by a periodic motion on an interval.

EXAMPLE 2.2. Suppose $T = \mathbb{R}$. The collection of maps

$$(2.6) \qquad \phi_t(z, v) = (z, \{v + s(z)t\}_{q(z)}), \qquad t \in \mathbb{R},$$

where $s(z) \in \mathbb{R} \setminus \{0\}$, $q(z) \in \mathbb{R}_+$ a.e. are measurable functions, is a measurable flow on $Z \times [0, q(\cdot))$. It is a cyclic flow because each point of the space comes back to its initial position in a finite (nonzero) time. We may think of the function $|s(z)|$ as the *speed* at which a point $(z, v)$ moves under the flow $\{\phi_t\}_{t \in \mathbb{R}}$. Observe also that the flow $\{\phi_t\}_{t \in \mathbb{R}}$ is isomorphic to the flow $\widetilde{\phi}_t(z, v) = (z, \{v + t\}_{q(z)/|s(z)|})$ on $Z \times [0, q(\cdot)/|s(\cdot)|)$.

EXAMPLE 2.3. Consider now the case $T = \mathbb{Z}$. The maps

$$(2.7) \qquad \phi_t(z, v) = (z, \{v + s(z)t\}_{q(z)}), \qquad t \in \mathbb{Z},$$

still define a cyclic flow on $Z \times ([0, q(\cdot)) \cap \mathbb{Z})$ or, equivalently, on the space $(Z \times [0, q(\cdot)), \sigma(dz)\lambda(dv))$ with $\lambda(dv) = \delta_{\mathbb{Z}}(dv)$ by using the notation of Definition 2.2. The definition of this cyclic flow, however, is not very natural. Consider, for example, the flow $\phi_t(v) = \{v + 2t\}_4$, $t \in \mathbb{Z}$, defined by (2.7) with the suppressed $Z = \{1\}, \sigma(dz) = \delta_{\{1\}}(dz)$ on the space $\{0, 1, 2, 3\}$. Since $t \in \mathbb{Z}$, this flow takes 0 to 2 and then 2 back to 0, and takes 1 to 3 and then 3 back to 1. It hence consists of two separate cyclic flows: the flow $\phi|_{\{0,2\}}$ restricted to the points $\{0, 2\}$ and the flow $\phi|_{\{1,3\}}$ restricted to the points $\{1, 3\}$. For a fixed $z$, the flow in the $v$ coordinate of (2.7) may hence consist of a number of distinct cyclic flows.

To avoid this type of situation when $T = \mathbb{Z}$, it is preferable to consider, instead of (2.7), the flow

$$(2.8) \qquad \phi_t(z, v) = (z, \{v + s(z)t\}_{|s(z)|q(z)}), \qquad t \in \mathbb{Z},$$

on $Z \times ([0, |s(\cdot)|q(\cdot)) \cap |s(\cdot)|\mathbb{Z})$, where $q(z) \in \mathbb{Z}_+$ a.e. and $a\mathbb{Z} = \{ap : p \in \mathbb{Z}\}$, $a \in \mathbb{R}$, or equivalently, on the space $(Z \times [0, |s(\cdot)|q(\cdot)), \sigma(dz)\delta_{|s(\cdot)|\mathbb{Z}}(dv))$. For example, the flow $\phi_t(v) = \{v + 2t\}_4$ is now defined only on the points $\{0, 2\}$ and thus, we do not have to deal with two flows anymore. In general, for a fixed $z$, the flow in the $v$ coordinate of (2.8) takes 0 to $|s(z)|$, $|s(z)|$ to $2|s(z)|, \ldots, (q(z) - 2)|s(z)|$ to $(q(z) - 1)|s(z)|$ before returning to 0. Since the space $[0, |s(z)|q(z)) \cap |s(z)|\mathbb{Z}$ consists only of these points $0, |s(z)|, \ldots, (q(z) - 1)|s(z)|$, there can be no other distinct cyclic flow on this space. Observe that the function $|s(z)|$ in (2.8) still plays the role of speed.

Observe also that the flows (2.6) and (2.8) have a common representation $(z, \{v + s(z)t\}_{q(z)})$ only for $|s(z)| = 1$. Since we prefer to work with a cyclic flow representation valid for both $T = \mathbb{R}$ and $T = \mathbb{Z}$, and since a flow $(z, \{v + s(z)t\}_{q(z)})$ with $|s(z)| = 1$ is isomorphic to the simpler cyclic flow $(z, \{v + t\}_{q(z)})$, we suppose in (2.4) of Definition 2.2 that $s(z) = 1$.

THEOREM 2.1. *Definitions 2.1 and 2.2 of cyclic flows are equivalent.*

PROOF. If the flow $\{\phi_t\}_{t \in T}$ is cyclic in the sense of Definition 2.2, then every point in the space $Z \times [0, q(\cdot))$ is cyclic and hence $S = L$ a.e. by using (2.5).

To show the converse, we suppose that $\{\phi_t\}_{t \in T}$ is cyclic in the sense of Definition 2.1. We first consider the case $T = \mathbb{R}$. Since $\{\phi_t\}_{t \in \mathbb{R}}$ has no fixed points (a.e.), we may suppose without loss of generality that the flow $\{\phi_t\}_{t \in \mathbb{R}}$ is a special flow on a space $Y \times [0, r(\cdot))$ as defined in Appendix A (see also Figure 1 in that appendix), that is,

$$(2.9) \qquad \phi_t(y, u) = (V^n y, t + u - r_n(y))$$

for $r_n(y) < u + t \leq r_{n+1}(y)$, where $r_n(y) = \sum_{k=0}^{n-1} r(V^k y)$, $n \geq 1$, $r_0(y) = 0$, $r_n(y) = \sum_{k=n}^{-1} r(V^k y)$, $n \leq -1$, $r(\cdot) > 0$ a.e. and $V$ is a one-to-one, onto, bimeasurable map on a a standard Lebesgue space $(Y, \tau)$. Indeed, as stated in Appendix A, given a flow $\{\phi_t\}_{t \in \mathbb{R}}$ without fixed points, there is a special flow given by (2.9) which is null isomorphic to $\{\phi_t\}_{t \in \mathbb{R}}$. If $\{\phi_t\}_{t \in \mathbb{R}}$ is cyclic in the sense of Definition 2.1, then the null isomorphic special flow (2.9) is cyclic as well. Then, if the flow (2.9) is shown to be null isomorphic to a flow given by (2.4), then the original flow $\{\phi_t\}_{t\in\mathbb{R}}$ is null isomorphic to the flow (2.4) as well (this is because null isomorphism is an equivalence relationship).

Since, by assumption, a.e. point $(y, u)$ comes back to its initial position in a finite period of time, we have that

$$
\begin{aligned}
Y &= \bigcup_{n=1}^{\infty} \{y : V^n y = y\} =: \bigcup_{n=1}^{\infty} A'_n \\
&= \bigcup_{n=1}^{\infty} (A'_n \setminus (A'_1 \cup \cdots \cup A'_{n-1})) =: \bigcup_{n=1}^{\infty} A_n \qquad (\text{with } A'_0 = \varnothing)
\end{aligned}
\tag{2.10}
$$

a.e. $\tau(dy)$. The set $A_n$ represents those $y$ that are attained for the first time by $V^n y$. Since $V A_n = A_n$, the sets $C_n := \{(y, u) : y \in A_n\}$ are invariant under the flow. We now want to show that the flow (2.9) satisfying (2.10) is indeed cyclic in the sense of Definition 2.2. Since the sets $C_n$ are invariant under the flow, it is enough to show that, for $n \geq 1$, the flow (2.9) restricted to $C_n$ is cyclic. We do so for $n = 2$ only, since the proof for other values of $n$ is similar.

To prove that $\phi|_{C_2}$ is cyclic, we show that there is a null isomorphism mapping $\phi|_{C_2}$ into a flow $\widetilde{\phi}$ of the type (2.4). The first step is to construct a space where the flow $\widetilde{\phi}$ is defined. The idea is to reduce the space so that $(y, v)$ and $(Vy, v)$ are only represented by either $(y, v)$ or $(Vy, v)$ or, since we are focusing on $A_2$, to reduce the space so that $y \in A_2$ and $Vy$ are represented only by either $y$ or $Vy$. To do so, we proceed in a way similar to the exhaustion principle used in ergodic theory [see, e.g., page 17 in Krengel (1985)]. Let $\widetilde{\tau}$ be a finite measure on $Y$ equivalent to $\tau$. Let first

$$\mathcal{B}_1 = \{\text{measurable } B \subset A_2 : B \cap VB = \varnothing\},$$

$\widetilde{\tau}_1 = \sup\{\widetilde{\tau}(B) : B \in \mathcal{B}_1\}$ and take $B_1 \in \mathcal{B}_1$ such that $\widetilde{\tau}(B_1) \geq \widetilde{\tau}_1/2$. Then define a sequence of sets $B_n$, $n \geq 2$, recursively, by letting

$$\mathcal{B}_n = \{\text{measurable } B \subset A_2 \setminus (B_1 \cup \cdots \cup B_{n-1} \cup VB_1 \cup \cdots \cup VB_{n-1}) : B \cap VB = \varnothing\},$$

$\widetilde{\tau}_n = \sup\{\widetilde{\tau}(B) : B \in \mathcal{B}_n\}$ and picking $B_n \in \mathcal{B}_n$ such that $\widetilde{\tau}(B_n) \geq \widetilde{\tau}_n/2$. Since $\widetilde{\tau}$ is finite and the sets $B_1, \ldots, B_{n+1}$ are disjoint, we have $\widetilde{\tau}(B_n) \to 0$ and hence $\widetilde{\tau}_n \to 0$ as $n \to \infty$.

We argue next that

$$A_2 = \bigcup_{n=1}^{\infty} (B_n \cup VB_n) \tag{2.11}$$

a.e. $\widetilde{\tau}(dy)$ and hence a.e. $\tau(dy)$.

Relationship (2.11) must hold because if it does not then we have a contradiction: we will show that it is then possible to construct a measurable set $B \subset A_2$ with $\widetilde{\tau}(B) > 0$, $B \cap VB = \varnothing$ and the sets $B$ and $B_n \cup VB_n$ being disjoint for all $n$. This is a contradiction because the argument preceding (2.11) precludes the existence of such a set. Assume then that (2.11) does not hold, that is, there is a set $A \subset A_2$ such that $\widetilde{\tau}(A) > 0$ and the sets $A$ and $\bigcup_{n=1}^{\infty}(B_n \cup VB_n)$ are disjoint. By the definition of $A_2$, we have $Vy \neq y$ on $A$ (a.e.) and hence there is a set $A_0 \subset A$ such that $\widetilde{\tau}(A_0) > 0$ and $A_0 \neq VA_0$ a.e. Then define $B = A_0 \setminus VA_0$. Since $A_0 \neq VA_0$ a.e., we have $\widetilde{\tau}(B) > 0$. Moreover, since $B = A_0 \cap VA_0^c$, we have $B \cap VB = A_0 \cap VA_0^c \cap VA_0 \cap V^2A_0^c = A_0 \cap V(A_0^c \cap A_0) \cap V^2A_0^c = \varnothing$ because $V$ is one-to-one. Since $A$ and $\bigcup_{n=1}^{\infty}(B_n \cup VB_n)$ are disjoint, the sets $B \subset A$ and $B_n \cup VB_n$ are disjoint for all $n$ as well. We therefore conclude that (2.11) holds.

Now let

$$Z = \bigcup_{n=1}^{\infty} B_n$$

and observe that

$$A_2 = \bigcup_{n=1}^{\infty}(B_n \cup VB_n) = \left(\bigcup_{n=1}^{\infty} B_n\right) \cup \left(\bigcup_{n=1}^{\infty} VB_n\right) = Z \cup VZ$$

a.e. $\tau(dy)$.

The spaces $Z$ and $VZ$ are disjoint by construction. Instead of focusing on the space $A_2 \subset Y$, where $A_2 = Z \cup VZ$ a.e. $\tau(dy)$, we "combine" $Z$ and $VZ$. We do so by focusing on the space $Z$ only and by defining the flow $\widetilde{\phi}$ on the space $(Z \times [0, r_2(\cdot)), \tau(dz)\,dv)$ as

(2.12) $$\widetilde{\phi}_t(z,v) = (z, \{v+t\}_{r_2(z)}),$$

where

$$r_2(z) = r(z) + r(Vz)$$

is the function which appears in the special representation (2.9). (To visualize this, see Figure 1 in Appendix A.) By using $r_2(z)$, we have replaced the "vertical" motions on $\{(z,u), 0 \leq u < r(z)\}$ and $\{(Vz,u), 0 \leq u < r(Vz)\}$ by a single motion on

$$\{(z,u), 0 \leq u < r(z) + r(Vz)\}.$$

Our new space $Z \times [0, r_2(\cdot))$ is thus related to the previous space $C_2 = \{(y,u) : y \in A_2\}$ by the map $\Phi : Z \times [0, r_2(\cdot)) \mapsto C_2$ defined by

$$\Phi(z,v) = \begin{cases} (z,v), & \text{if } 0 \leq v < r(z), \\ (Vz, v - r(z)), & \text{if } r(z) \leq v < r_2(z). \end{cases}$$

Then, by using $A_2 = Z \cup VZ$ a.e. $\tau(dy)$ and the fact that $Z$ and $VZ$ are disjoint, we obtain that $\Phi$ is a null isomorphism and

$$\phi_t|_{C_2}(\Phi(z,v)) = \Phi(\widetilde{\phi}_t(z,v)).$$

This shows that the flow $\phi|_{C_2}$ is indeed cyclic in the sense of Definition 2.2.

To show the converse in the case $T = \mathbb{Z}$ is easier. We sketch here only the main ideas of the proof. The function $q(s) = \min\{n \in \mathbb{N} : \phi_n(s) = s\}$ is measurable and a.e. finite on $S$. It is enough to show, for example, that the flow $\{\phi_t\}_{t \in \mathbb{Z}}$ restricted to $C_2 = \{s : q(s) = 2\}$ is cyclic. Arguing as above, we can construct a measurable set $Z$ such that $C_2 = Z \cup \phi_1 Z$, where $Z$ and $\phi_1 Z$ are disjoint. The flow $\phi|_{C_2}$ can then be shown to be isomorphic to the flow $\widetilde{\phi}_t(z,v) = (z, \{v+t\}_2)$. $\square$

In the following lemma, we characterize cocycles associated with cyclic flows. (See the end of Appendix A for a definition of a cocycle.) This result is used in the next section.

LEMMA 2.2. *Let $\{\phi_t\}_{t \in T}$ be a cyclic flow and let $\{a_t\}_{t \in T}$ be a cocycle for $\{\phi_t\}_{t \in T}$ taking values in a second countable group $(A, \cdot)$. Suppose that $\Phi : Z \times [0, q(\cdot)) \setminus \widetilde{N} \mapsto S \setminus N$ is the null isomorphism between the flows $\{\phi_t\}_{t \in T}$ and $\{\widetilde{\phi}_t\}_{t \in T}$ in Definition 2.2. Let $\widetilde{a}_t(z,v) = a_t(\Phi(z,v))$ if $(z,v) \in Z \times [0, q(\cdot)) \setminus \widetilde{N}$, and $\widetilde{a}_t(z,v) = e$ if $(z,v) \in \widetilde{N}$, where $e$ is the group unity. Then $\{\widetilde{a}_t\}_{t \in T}$ is a cocycle for $\{\widetilde{\phi}_t\}_{t \in T}$ and*

$$\widetilde{a}_t(z,v) = (\widetilde{a}(z,v))^{-1} \widetilde{a}_1(z)^{[v+t]_{q(z)}} \widetilde{a}(z, \{v+t\}_{q(z)}) \tag{2.13}$$

*for all $t \in T$ and $(z,v) \in Z \times [0, q(\cdot))$, where $\widetilde{a} : Z \times [0, q(\cdot)) \mapsto A$ and $\widetilde{a}_1 : Z \mapsto A$ are some measurable functions.*

PROOF. We may suppose without loss of generality that $N = \widetilde{N} = \varnothing$ because the proof below shows that $\widetilde{a}_t(z,v) = a_t(\Phi(z,v))$ satisfies the cocycle relationship on the set $Z \times [0, q(\cdot)) \setminus \widetilde{N}$ (which is invariant for the flow $\widetilde{\phi}$) and so obviously does $\widetilde{a}_t(z,v) = e$ on the set $\widetilde{N}$. By substituting $s = \Phi(z,v)$ in the definition (A.3) of a cocycle, we obtain that

$$a_{t_1+t_2}(\Phi(z,v)) = a_{t_2}(\Phi(z,v)) a_{t_1}(\phi_{t_2}(\Phi(z,v)))$$

and hence, since $\phi_t \circ \Phi = \Phi \circ \widetilde{\phi}_t$, we get that

$$\widetilde{a}_{t_1+t_2}(z,v) = \widetilde{a}_{t_2}(z,v) \widetilde{a}_{t_1}(\widetilde{\phi}_{t_2}(z,v)) \tag{2.14}$$

$$= \widetilde{a}_{t_2}(z,v) \widetilde{a}_{t_1}(z, \{v+t_2\}_{q(z)}). \tag{2.15}$$

Relationship (2.14) shows that $\{\widetilde{a}_t\}_{t \in T}$ is a cocycle for the flow $\{\widetilde{\phi}_t\}_{t \in T}$. To show (2.13), we use (2.15). We consider the case $Z = \{1\}$ only. The proof for a general space $Z$ follows as below by working with a fixed $z$. For notational simplicity, we denote $\widetilde{a}_t(1,v)$ by $\widetilde{a}_t(v)$ and, to avoid writing indices, by $\widetilde{a}(t,v)$.

By taking $v = 0$ in (2.15), we get $\widetilde{a}(t_1 + t_2, 0) = \widetilde{a}(t_2, 0)\widetilde{a}(t_1, \{t_2\}_q)$. Then

$$(2.16) \qquad \widetilde{a}(t,v) = (\widetilde{a}(v,0))^{-1}\widetilde{a}(v+t,0)$$

if $v \in [0,q) \cap T$. Observe now that, by (1.16) and (2.15),

$$\widetilde{a}(v+t,0) = \widetilde{a}(q[v+t]_q + \{v+t\}_q, 0) = \widetilde{a}(q[v+t]_q, 0)\widetilde{a}(\{v+t\}_q, 0)$$

for all $t \in T$ and $v \in [0,q) \cap T$. Then, by (2.16), for $v \in [0,q) \cap T$,

$$\widetilde{a}(t,v) = (\widetilde{a}(v,0))^{-1}\widetilde{a}(q[v+t]_q, 0)\widetilde{a}(\{v+t\}_q, 0)$$
$$= (\widetilde{a}(v))^{-1}\widetilde{a}(q[v+t]_q)\widetilde{a}(\{v+t\}_q),$$

where $\widetilde{a}(\cdot) = \widetilde{a}(\cdot,0)$, but by setting $t_2 = nq$, $t_1 = mq$ and $v = 0$ in (2.15), we get that $\widetilde{a}(qm + qn) = \widetilde{a}(qm)\widetilde{a}(qn)$ for all $n, m \in \mathbb{Z}$. It follows that $\widetilde{a}(qm) = \widetilde{a}_1^m$ for some $\widetilde{a}_1 \in A$ and hence

$$\widetilde{a}(t,v) = (\widetilde{a}(v))^{-1}\widetilde{a}_1^{[v+t]_q}\widetilde{a}(\{v+t\}_q),$$

which proves (2.13) when $Z = \{1\}$. □

**3. Representation of stationary processes generated by periodic flows.** We now provide a representation of stable stationary processes generated by periodic flows. This basic result is used several times in this and the following section.

THEOREM 3.1. *Suppose that a stationary $S\alpha S$, $\alpha \in (0,2)$, process $X_\alpha$ is generated by a periodic flow in the sense of Definition 1.1. Then $X_\alpha$ can be represented in distribution as the sum of two independent stationary stable processes. The first process is a harmonizable process in the complex-valued case or a trivial process in the real-valued case,*

$$(3.1) \qquad \begin{aligned} &\int_{\widehat{T}} e^{itx} N_\alpha(dx) &&\text{(complex-valued case),} \\ &\int_{\{1\}} N_\alpha(dx) + n(t)\int_{\{2\}} N_\alpha(dx) &&\text{(real-valued case),} \end{aligned}$$

*where $N_\alpha$ has a finite control measure $\eta(dx)$ and $n(t)$ is defined by (1.12). The second process can be represented as*

$$(3.2) \qquad \int_Z \int_{[0,q(z))} b_1(z)^{[v+t]_{q(z)}} g(z, \{v+t\}_{q(z)}) M_\alpha(dz, dv).$$

*Here, $(Z, \mathcal{Z}, \sigma)$ is a standard Lebesgue space, $q(z) \in T_+$ a.e. $\sigma(dz)$, $g \in L^\alpha(Z \times [0, q(\cdot)), \sigma(dz)\lambda(dv))$ and $b_1(z) \in \{-1,1\}$ [or $b_1(z) \in \{w : |w| = 1\}$ in the complex-valued case] are measurable functions, and $S\alpha S$ random measure $M_\alpha$ has the control measure $\sigma(dz)\lambda(dv)$.*

PROOF. Suppose that the process $X_\alpha$ is generated by a flow $\{\phi_t\}_{t\in T}$ on $(S,\mathcal{S},\nu)$ which is periodic. Since the flow is periodic, we have $S = P$ a.e. $\nu(ds)$ and hence $S = F + L$ a.e. $\nu(ds)$ as well, where $F$ and $L$ are the fixed and the cyclic points of the flow $\{\phi_t\}_{t\in T}$. Then

$$X_\alpha(t) \stackrel{d}{=} \int_F f_t(s) M_\alpha(ds) + \int_L f_t(s) M_\alpha(ds) =: Y_\alpha(t) + Z_\alpha(t),$$

where the stationary stable processes $Y_\alpha$ and $Z_\alpha$ are independent, the process $Y_\alpha$ is generated by an identity flow and the process $Z_\alpha$ is generated by a cyclic flow. (The processes $Y_\alpha$ and $Z_\alpha$ are independent because $F \cap L = \varnothing$.) By Propositions 5.1 and 5.2 in Rosiński (1995), the process $Y_\alpha$ is harmonizable (or trivial). To conclude the theorem, we still need to show that the process $Z_\alpha$ has a representation (3.2).

By Definition 2.2, there is a space $(Z, \mathcal{Z}, \sigma)$, function $q(z) \in T_+$ a.e. $\sigma(dz)$ and a null isomorphism $\Phi : Z \times [0, q(\cdot)) \mapsto L$ such that

$$(3.3) \qquad \phi_t(\Phi(z, v)) = \Phi(z, \{v + t\}_{q(z)})$$

for all $t \in T$ and $(z, v) \in Z \times [0, q(\cdot))$. In other words, the flow $\{\phi_t\}_{t \in T}$ on $(L, \nu)$ is null isomorphic to the flow $\{\widetilde{\phi}_t\}_{t \in T}$ on $(Z \times [0, q(\cdot)), \sigma(dz)\lambda(dv))$ defined by

$$\widetilde{\phi}_t(z, v) = (z, \{v + t\}_{q(z)}).$$

[We may suppose that the null sets in (2.5) are empty because, otherwise, we can replace $L$ by $L \setminus N$ in the definition of $Z_\alpha$ without changing its distribution.] By replacing $s$ with $\Phi(z, v)$ in (1.4) and using (3.3), we get that for all $t \in T$,

$$(3.4) \quad f_t(\Phi(z,v)) = a_t(\Phi(z,v)) \left\{ \frac{d(\nu \circ \phi_t)}{d\nu}(\Phi(z,v)) \right\}^{1/\alpha} f_0(\Phi(\widetilde{\phi}_t(z,v)))$$

a.e. $\sigma(dz)\lambda(dv)$. Now, by Lemma 2.2, $a_t(\Phi(z,v)) = \widetilde{a}_1(z)^{[v+t]_{q(z)}} \widetilde{a}(\widetilde{\phi}_t(z,v))/\widetilde{a}(z,v)$. Since $\phi_t \circ \Phi = \Phi \circ \widetilde{\phi}_t$, we also have that

$$\frac{d(\nu \circ \phi_t)}{d\nu} \circ \Phi = \frac{d(\nu \circ \Phi \circ \widetilde{\phi}_t)}{d(\nu \circ \Phi)}$$

$$= \left( \frac{d\nu}{d((\sigma \otimes \lambda) \circ \Phi^{-1})} \circ \Phi \circ \widetilde{\phi}_t \right) \frac{d((\sigma \otimes \lambda) \circ \widetilde{\phi}_t)}{d(\sigma \otimes \lambda)} \frac{d(\sigma \otimes \lambda)}{d(\nu \circ \Phi)}$$

$$= \left( \frac{d\nu}{d((\sigma \otimes \lambda) \circ \Phi^{-1})} \circ \Phi \circ \widetilde{\phi}_t \right) \frac{d((\sigma \otimes \lambda) \circ \Phi^{-1})}{d\nu} \circ \Phi$$

$$= \left( \frac{d\nu}{d((\sigma \otimes \lambda) \circ \Phi^{-1})} \circ \Phi \circ \widetilde{\phi}_t \right) \left( \frac{d\nu}{d((\sigma \otimes \lambda) \circ \Phi^{-1})} \circ \Phi \right)^{-1},$$

where $d((\sigma \otimes \lambda) \circ \widetilde{\phi}_t)/d(\sigma \otimes \lambda) \equiv 1$ because the first component in $\widetilde{\phi}_t(z,v) = (z, \{v+t\}_{q(z)})$ remains the same and the second, where $v$ is the variable, preserves the measure $\lambda$. Hence, by setting

$$(3.5) \quad g_t(z,v) = \widetilde{a}(z,v)\left\{\frac{d\nu}{d((\sigma \otimes \lambda) \circ \Phi^{-1})}(\Phi(z,v))\right\}^{1/\alpha} f_t(\Phi(z,v))$$

in relationship (3.4), we obtain that, for all $t \in T$,

$$(3.6) \quad g_t(z,v) = \widetilde{a}_1(z)^{[v+t]_{q(z)}} g_0(\widetilde{\phi}_t(z,v))$$

a.e. $\sigma(dz)\lambda(dv)$. Finally, observe that by writing the characteristic functions, it is easy to see that (3.5) implies

$$\{Z_\alpha(t)\}_{t \in T} \stackrel{d}{=} \left\{\int_L f_t(s) M_\alpha(ds)\right\}_{t \in T} \stackrel{d}{=} \left\{\int_Z \int_{[0,q(z))} g_t(z,v) \widetilde{M}_\alpha(dz,dv)\right\}_{t \in T},$$

where $\widetilde{M}_\alpha(dz,dv)$ has the control measure $\sigma(dz)\lambda(dv)$. The result of the theorem then follows from (3.6) by setting $b_1(z) = \widetilde{a}_1(z)$ and $g(z,v) = g_0(z,v)$. □

REMARK 3.1. The proof of Theorem 3.1 shows that stationary S$\alpha$S processes generated by cyclic flows have a representation (3.2).

DEFINITION 3.1. A stationary stable process that can be represented by the sum of the processes (3.1) and (3.2) as in Theorem 3.1 is called *a stationary periodic process*.

The following result is useful for recognizing stationary periodic processes and, more specifically, processes (3.2) when $T = \mathbb{R}$.

PROPOSITION 3.1. *With the notation of Theorem 3.1 and letting $s(z) \in \mathbb{R} \setminus \{0\}$ a.e. be a measurable function, processes*

$$(3.7) \quad \int_Z \int_{[0,q(z))} b_1(z)^{[v+s(z)t]_{q(z)}} g(z, \{v+s(z)t\}_{q(z)}) M_\alpha(dz,dv), \qquad t \in \mathbb{R},$$

*have a representation (3.2) (with possibly different functions $q$ and $g$) and hence are stationary periodic processes.*

PROOF. By using the relationships

$$\{v+st\}_q = q - \{-v-st\}_q = q - \{(q-v) - st\}_q,$$
$$[v+st]_q = \frac{1}{q}(v+st - \{v+st\}_q)$$
$$= -\frac{1}{q}((q-v) - st - \{(q-v) - st\}_q) = -[(q-v) - st]_q$$

and by making the change of variables $v$ to $q(z) - v$ when $s(z) < 0$, we can first represent the process in (3.7) as

$$(3.8) \quad \int_Z \int_{[0,q(z))} b_1(z)^{[v+|s(z)|t]_{q(z)}} \widehat{g}(z, \{v + |s(z)|t\}_{q(z)}) M_\alpha(dz, dv),$$

where $\widehat{g}(z, u) = g(z, q(z) - u)$ if $s(z) < 0$ and $\widehat{g}(z, u) = g(z, u)$ if $s(z) > 0$. Then, by using the relationships

$$\{v + |s|t\}_q = |s|\{|s|^{-1}v + t\}_{|s|^{-1}q},$$

$$[v + |s|t]_q = \frac{1}{q}((v + |s|t) - \{v + |s|t\}_q)$$

$$= \frac{1}{|s|^{-1}q}((|s|^{-1}v + t) - \{|s|^{-1}v + t\}_{|s|^{-1}q}) = [|s|^{-1}v + t]_{|s|^{-1}q}$$

and by making the change of variables $|s(z)|^{-1}v = \widetilde{v}$, we can represent the process in (3.8) as

$$(3.9) \quad \int_Z \int_{[0,\widetilde{q}(z))} b_1(z)^{[\widetilde{v}+t]_{\widetilde{q}(z)}} \widetilde{g}(z, \{\widetilde{v} + t\}_{\widetilde{q}(z)}) M_\alpha(dz, d\widetilde{v}),$$

where $\widetilde{g}(z, u) = |s(z)|^{1/\alpha} g(z, |s(z)|^{-1}u)$ and $\widetilde{q}(z) = |s(z)|^{-1} q(z)$.  □

REMARK 3.2. We can also use Theorem 3.1 to show that the process (3.7) has a representation (3.2). Indeed, Example 2.2 shows that $\phi_t(z, v) = (z, \{v + s(z)t\}_{q(z)})$, $t \in \mathbb{R}$, is a cyclic flow on $Z \times [0, q(\cdot))$. By using the relationship $[v + s(t_1 + t_2)]_q = [v + st_1]_q + [\{v + st_1\}_q + st_2]_q$ [to verify it, use the second relationship in (1.16) and the fact that $\{\phi_t\}$ is a flow], we get that $b_1(z)^{[v+s(z)t]_{q(z)}}$ is a cocycle for the flow $\{\phi_t\}_{t \in \mathbb{R}}$. The process (3.7) is thus generated by the flow $\{\phi_t\}_{t \in \mathbb{R}}$ in the sense of Definition 1.1. Since this flow is cyclic, the remark before Definition 3.1 shows that the process (3.7) has a representation (3.2). Observe, however, that this does not prove that to obtain (3.2), only the functions $q$ and $g$ may need to be modified.

The term "periodic" in "stationary periodic processes" refers to a process that has a representations (3.1) and (3.2), where the kernel has a periodic-like structure [as in (4.1)]. It does not necessarily imply that an underlying generating flow of the process is periodic. In fact, as the following elementary example shows, without any restrictions on the kernel of a process, a stationary periodic process can be generated in the sense of Definition 1.1 by conservative flows other than periodic flows.

EXAMPLE 3.1. Let $(Y, \mathcal{Y}, \tau)$ be a standard Lebesgue space with $0 < \tau(Y) < \infty$. Observe that a stationary periodic process that has a representation (3.2) can also be represented as

$$(3.10) \quad (\tau(Y))^{-1/\alpha} \int_Y \int_Z \int_{[0,q(z))} b_1(z)^{[v+t]_{q(z)}} g(z, \{v + t\}_{q(z)}) M_\alpha(dy, dz, dv),$$

where the S$\alpha$S random measure $M_\alpha$ has control measure $\tau(dy)\sigma(dz)\lambda(dv)$, because there is no variable $y$ in the kernel of (3.10). Let now $\{\phi_t\}_{t\in\mathbb{R}}$ be any measure preserving conservative flow on $(Y, \mathcal{Y}, \tau)$. Then the stationary periodic process (3.2), when represented by (3.10), is also generated by the flow $\widetilde{\phi}_t(y, z, v) = (\phi_t(y), z, \{v+t\}_{q(z)})$ on $Y \times Z \times [0, q(\cdot))$ in the sense of Definition 1.1. The generating flow is therefore not unique. Observe that since we can choose the flow $\{\phi_t\}_{t\in\mathbb{R}}$ to be nonperiodic, the flow $\{\widetilde{\phi}_t\}_{t\in\mathbb{R}}$ will also be nonperiodic. A similar problem exists when we consider harmonizable (or trivial) processes and identity flows.

Without any restrictions on a kernel function, the generating flow may not be unique. In this case, not only stationary periodic processes can be generated by nonperiodic flows, but harmonizable (or trivial) processes can also be represented by (3.2) (Lemma 3.1). This result further indicates that we cannot associate, in general, harmonizable (or trivial) processes with identity flows and processes that have the representation (3.2) with cyclic flows.

LEMMA 3.1. *The S$\alpha$S, $\alpha \in (0, 2)$, harmonizable processes (or trivial processes in the real-valued case) can be represented as* (3.2).

PROOF. Consider the process

$$
\begin{aligned}
X_\alpha(t) &= \int_{\widehat{T}} \int_{[0,2)} (e^{i2z})^{[v+t]_2} (e^{iz})^{\{v+t\}_2} M_\alpha(dz, dv) \\
&= \int_{\widehat{T}} \int_{[0,2)} (e^{iz})^{(v+t)} M_\alpha(dz, dv), \qquad t \in T,
\end{aligned}
\tag{3.11}
$$

where $M_\alpha$ is a complex-valued rotationally invariant S$\alpha$S random measure with control measure $\eta(dz)\lambda(dv)$ and $\eta(\widehat{T}) < \infty$ [see also the notation (1.10)]. The process $X_\alpha$ has a representation (3.2) with $Z = \widehat{T}$, $\sigma(dz) = \eta(dz)$, $b_1(z) \equiv e^{i2z}$, $q(z) = 2$ and $g(z, u) = e^{izu}$. Observe that since $e^{iz(v+t)} = e^{izv}e^{izt}$, $|e^{izv}| = 1$ and $e^{izv}$ does not involve time $t$,

$$
(3.12)\ \{X_\alpha(t)\}_{t\in T} \stackrel{d}{=} \left\{ \int_{\widehat{T}} \int_{[0,2)} e^{izt} M_\alpha(dz, dv) \right\}_{t\in T} \stackrel{d}{=} \left\{ 2^{1/\alpha} \int_{\widehat{T}} e^{izt} M_\alpha(dz) \right\}_{t\in T},
$$

where $M_\alpha$ is a complex-valued rotationally invariant measure with control measure $\eta(dz)$. Hence, $X_\alpha$ is also a harmonizable process, showing the result for harmonizable processes.

The case of trivial processes with $T = \mathbb{R}$ follows by taking, for example, $Z = \{1\}$, $b_1(z) = 1$, $g(1, z) \equiv 1$ and $q(1) = 1$ in (3.2). When $T = \mathbb{Z}$, take $Z =$ $\{1,2\}$, $b_1(z) \equiv 1$ and $g(z,v) \equiv a(z)^v$ with $a(1) = 1, a(2) = -1$ and $q(z) = 2$. Then (3.2) becomes

$$\int_{\{1,2\}} \int_{\{0,1\}} (a(z))^{\{v+t\}_2} M_\alpha(dz, dv)$$

$$= \int_{\{1,2\}} \int_{\{0,1\}} (a(z))^{v+t} M_\alpha(dz, dv)$$

$$\stackrel{d}{=} \int_{\{1,2\}} a(z)^t N_\alpha(dz)$$

$$= \int_{\{1\}} N_\alpha(dz) + (-1)^t \int_{\{2\}} N_\alpha(dz),$$

which shows the result in the case of trivial processes when $T = \mathbb{Z}$. □

Lemma 3.1 has the following implication:

COROLLARY 3.1. *Stationary periodic processes can also be represented by* (3.2).

The representation of the process in Example 3.1 and the representation (3.2) of a harmonizable (or trivial) process in Lemma 3.1 have a built-in redundancy [e.g., there is no variable $y$ in the kernel of (3.10)]. If we eliminate redundancy and focus on *minimal representations* only, then by Theorem 3.2, stationary periodic process defined by (3.1) and (3.2) can only be generated by periodic flows. This explains our use of the term "stationary periodic processes." (Another justification is provided in the following sections.) Moreover, as shown in the following theorem, under minimal representations, harmonizable (or trivial) processes cannot have a minimal representation (3.2) and they are generated only by identity flows.

THEOREM 3.2. *If representations* (3.1) *and* (3.2) *of* $X_\alpha$ *is minimal, then* $X_\alpha$ *is generated by a unique flow in the sense of Definition* 1.1. *This flow is periodic for* (3.1) *and* (3.2), *identity for* (3.1) *and cyclic for* (3.2). *The representation* (3.1) *is always minimal.*

PROOF. If the representation is minimal, then the generating flow is unique by Theorem 3.1 in Rosiński (1995). Representations (3.1) and (3.2) are obviously generated by an identity flow and a cyclic flow, respectively, in the sense of Definition 1.1 and, therefore, representations (3.1) and (3.2) of their sum is generated by a periodic flow in the sense of Definition 1.1. The minimality of the representation (3.1) in the complex-valued case was shown by Rosiński (1998a), Example 4.8. The minimality in the real-valued case can be seen directly from the definition of minimal representations. [In the case $T = \mathbb{R}$, since $n(t) \equiv 0$, we assume implicitly that the representation (3.1) is defined on the space $\{1\}$ and not $\{1,2\}$.] $\square$

By Definition 3.1, harmonizable (or trivial) processes are also stationary periodic processes. Here is another example of stationary stable process which is a stationary periodic process.

EXAMPLE 3.2. Consider the process

$$X_\alpha(t) = \int_\mathbb{R} \int_0^{2\pi} \cos(v + zt) M_\alpha(dz, dv), \qquad t \in \mathbb{R},$$

where the S$\alpha$S random measure $M_\alpha$ has control measure $\mu(dz)\,dv$ and $\mu(\mathbb{R}) < \infty$. The process $X_\alpha$ is well defined, that is, $\cos(v+zt) \in L^\alpha(\mathbb{R} \times (0, 2\pi), \mu(dz)\,dv)$ for each $t \in \mathbb{R}$. Since $\cos(u) = \cos(\{u\}_{2\pi})$, it has a representation (3.7) with $Z = \mathbb{R}$, $\sigma(dz) = \mu(dz)$, $b_1(z) \equiv 1$, $s(z) = z$, $q(z) = 2\pi$ and $g(z, u) = \cos(u)$. Hence, by Proposition 3.1, the process $X_\alpha$ is a stationary periodic process. This can also be seen directly by using the proof of Proposition 3.1 to observe that

$$X_\alpha(t) \stackrel{d}{=} \int_\mathbb{R} \int_0^{2\pi/|z|} \cos\left(z\{w + t\}_{2\pi/|z|}\right) M_\alpha(dz, dv).$$

As shown in Example 2.5 of Rosiński (2000), the process $X_\alpha$ has the same (up to a constant) finite-dimensional distributions as the real part of a harmonizable process (1.9); more precisely,

(3.13) $$\{X_\alpha(t)\}_{t\in\mathbb{R}} \stackrel{d}{=} \left\{(2\pi)^{1/\alpha} \Re \int_\mathbb{R} e^{itx} M_\alpha(dx)\right\}_{t\in\mathbb{R}},$$

where $M_\alpha$ is a complex-valued rotationally invariant measure with the control measure $\mu(dx)$.

Since harmonizable (or trivial) processes are also stationary periodic processes, we may want to single out stationary periodic processes that do not have a harmonizable (or trivial) component. The following definition makes this precise. The introduced terminology is often used in the sequel, along with that of stationary periodic processes.

DEFINITION 3.2. A stationary periodic process is called a *stationary cyclic process* if it does not have a harmonizable (or trivial) component, that is, it cannot be represented as a sum of two independent stationary processes, one of which is a nondegenerate harmonizable (or trivial) process.

REMARK 3.3. A stationary periodic (harmonizable or trivial, resp.) process is defined as a process which can have representations (3.1) and (3.2) [(3.1), resp.]. A stationary cyclic process, however, cannot be defined as a process which can have the representation (3.2), because, by Lemma 3.1, a harmonizable (or trivial) process can also be represented as (3.2). It is necessary, therefore, to exclude explicitly the harmonizable (or trivial) component in Definition 3.2.

How can one determine whether a given stationary process is a stationary periodic or cyclic process? Example 3.1 and Lemma 3.1 show that, in general, it is not enough to examine whether an underlying *flow* of the process is periodic or cyclic. We can, however, identify these processes through underlying flows if their representations are minimal (see Theorem 3.2 and also Corollary 6.2).

Since minimal representations are typically not easy to determine in practice, we would like to have an identification criterion which does not rely on minimal representations. We can do so through periodic and cyclic component sets which we define next. We work now with the kernel of a stationary process itself rather than with a generating flow. Flows, however, are still used as a tool to obtain an identification result (see Theorem 4.1 and its proof). The identification results are used in Section 6 to establish a further decomposition of stationary stable processes.

**4. Characterization of stationary periodic processes.** Consider a stationary process with representation (1.1) involving the kernel $f_t$. We first provide a criterion on $f_t$ for the process to be a stationary periodic process.

DEFINITION 4.1. *A periodic component set* of a stationary stable process $X_\alpha$ that has a representation (1.1) is defined as

$$(4.1) \quad C_P = \{s \in S : \exists h = h(s) \in T \setminus \{0\} : f_{t+h}(s) = a(h,s) f_t(s)$$
$$\text{a.e. } \lambda(dt) \text{ for some } a(h,s) \neq 0\}.$$

Relationship (4.1) expresses physically the fact that, starting at $s \in C_P$, we come back to $f_t(s)$ after some time $h(s)$.

LEMMA 4.1. *A periodic component set $C_P$ in* (4.1) *is $\nu$-measurable. Moreover, the functions $h(s)$ and $a(s) = a(h(s), s)$ in* (4.1) *can be taken to be $\nu$-measurable as well.*

PROOF. We first show that the set $C_P$ is $\nu$-measurable. Observe that the condition (4.1) says that the ratio $f_{t+h}(s)/f_t(s)$ does not depend on $t$ and hence $C_P$ can be also expressed as

$$C_P = \{s \in S : \exists h = h(s) \neq 0 : f_{t_1+h}(s) f_{t_2}(s) = f_{t_2+h}(s) f_{t_1}(s) \text{ a.e. } \lambda(dt_1)\lambda(dt_2)\}.$$

To deal with the potential measurability problem raised by $\exists h$, consider the set

$$A = \{(s,h) \neq (s,0) : f_{t_1+h}(s)f_{t_2}(s) = f_{t_2+h}(s)f_{t_1}(s) \text{ a.e. } \lambda(dt_1)\lambda(dt_2)\}.$$

Since $A = \{(s,h) \neq (s,0) : F(s,h) = 0\}$, where the function

$$F(s,h) = \int_T \int_T \mathbb{1}_{\{f_{t_1+h}(s)f_{t_2}(s) \neq f_{t_2+h}(s)f_{t_1}(s)\}}(s,h,t_1,t_2)\lambda(dt_1)\lambda(dt_2)$$

is Borel measurable by the Fubini's theorem [use Theorem A in Halmos (1950), page 143, and the fact that the function $(t,s) \mapsto f_t(s)$ is Borel], we obtain that the set $A$ is Borel measurable. We can verify now that $C_P$ is the projection of the set $A$ on $s$, namely, that $C_P = \text{proj}_S A := \{s : \exists h : (s,h) \in A\}$. By using Lemma 4.2, the set $C_P$ is $\nu$-measurable and we can choose the function $h(s)$ in (4.1) to be $\nu$-measurable. The $\nu$-measurability of $a(s)$ follows since, for $s \in C_P$, $f_{t_1+h(s)}(s)f_{t_2}(s) = f_{t_2+h(s)}(s)f_{t_1}(s)$ a.e. $\lambda(dt_1)\lambda(dt_2)$ and hence $a(s) = f_{t_2+h(s)}(s)(f_{t_2}(s))^{-1}$ a.e. $\nu(ds)$ for some $t_2 \in T$. □

The following result characterizes stationary periodic processes.

THEOREM 4.1. *A $S\alpha S$, $\alpha \in (0,2)$, stationary process $X_\alpha$ given by* (1.1) *with* $\text{supp}\{f_t, t \in T\} = S$ *a.e.* $\nu(ds)$ *is a stationary periodic process if and only if*

$$C_P = S \qquad \text{a.e. } \nu(ds),$$

*where $C_P$ is the periodic component set defined in* (4.1).

PROOF. Suppose first that $X_\alpha$ is a S$\alpha$S process given by (1.1) and that it is a stationary periodic process. Then, by Corollary 3.1 following Lemma 3.1, $X_\alpha$ has a representation (3.2) on a space $Z \times [0, q(\cdot))$ and with a kernel function

$$g_t(z,v) = b_1(z)^{[v+t]_{q(z)}} g(z, \{v+t\}_{q(z)}).$$

Since $\{v + (t+q(z))\}_{q(z)} = \{v+t\}_{q(z)}$ and $[v+(t+q(z))]_{q(z)} = [v+t]_{q(z)} + 1$, we have

$$g_{t+h(z,v)}(z,v) = a(z,v)g_t(z,v)$$

for all $(z,v) \in Z \times [0, q(\cdot))$, where $h(z,v) = q(z)$ and $a(z,v) = b_1(z)$. By Theorem 1.1 in Rosiński (1995), there are measurable maps $\Phi : S \mapsto Z \times [0, q(\cdot))$ and $k : S \mapsto \mathbb{R} \setminus \{0\}$ (or $\mathbb{C} \setminus \{0\}$) such that, for a.e. $\nu(ds)$, $f_t(s) = k(s)g_t(\Phi(s))$ a.e. $\lambda(dt)$. Then, for a.e. $\nu(ds)$,

$$f_{t+h(\Phi(s))}(s) = k(s)g_{t+h(\Phi(s))}(\Phi(s)) = k(s)a(\Phi(s))g_t(\Phi(s)) = a(\Phi(s))f_t(s)$$

a.e. $\lambda(dt)$. This shows that $C_P = S$ a.e. $\nu(ds)$.

Suppose now that $X_\alpha$ is a stationary stable process given by (1.1) with $\operatorname{supp}\{f_t, t \in T\} = S$ a.e. $\nu(ds)$ and such that $C_P = S$ a.e. $\nu(ds)$. We want to show that $X_\alpha$ is a stationary periodic process. The proof involves a number of steps.

STEP 1. First, we show that one may suppose without loss of generality that the representation $\{f_t\}_{t \in T}$ of the process $X_\alpha$ is minimal with $C_P = S$ a.e. (minimal representations are defined in Appendix B). By Theorem 2.2(a), in Rosiński (1995) [due to Hardin (1982), Theorem 5.1], the process $X_\alpha$ has a minimal integral representation

$$\int_{\widetilde{S}} \widetilde{f}_t(\tilde{s}) \widetilde{M}_\alpha(d\tilde{s}), \tag{4.2}$$

where $(\widetilde{S}, \widetilde{\mathcal{S}}, \widetilde{\nu})$ is some standard Lebesgue space, $\{\widetilde{f}_t\}_{t \in T} \subset L^\alpha(\widetilde{S}, \widetilde{\mathcal{S}}, \widetilde{\nu})$ and $\widetilde{M}_\alpha$ has the control measure $\widetilde{\nu}$. Let $\widetilde{C}_P$ be the periodic component set of $X_\alpha$ defined through the representation $\{\widetilde{f}_t\}_{t \in T}$. To conclude the first step, it is enough to show that $\widetilde{C}_P = \widetilde{S}$ a.e. $\widetilde{\nu}(d\tilde{s})$. By Remark 2.5 in Rosiński (1995), there are Borel measurable maps $\Phi: S \mapsto \widetilde{S}$ and $k: S \mapsto \mathbb{R} \setminus \{0\}$ (or $\mathbb{C} \setminus \{0\}$) such that, for any $t \in T$,

$$f_t(s) = k(s) \widetilde{f}_t(\Phi(s)) \tag{4.3}$$

a.e. $\nu(ds)$ and

$$\widetilde{\nu} = \nu_k \circ \Phi^{-1}, \tag{4.4}$$

where $\nu_k(ds) = |k(s)|^\alpha \nu(ds)$. Since, for $s \in C_P$, $f_{t+h(s)}(s) = a(s) f_t(s)$ a.e. $\lambda(dt)$, we expect, in view of (4.3) that, for a.e. $s \in C_P$,

$$\widetilde{f}_{t+h(s)}(\Phi(s)) = a(s) \widetilde{f}_t(\Phi(s)) \qquad \text{a.e. } \lambda(dt) \tag{4.5}$$

and hence that $\Phi(s) \in \widetilde{C}_P$ a.e. $\nu(ds)$.

To demonstrate that (4.5) follows from (4.3), consider first the set $A = \{(s,h): f_{t+h}(s) = k(s) \widetilde{f}_{t+h}(\Phi(s))$ and $f_{t+h}(s) = a(h,s) f_t(s)$ a.e. $\lambda(dt)$ for some $a(h,s) \neq 0\}$ and, by Lemma 4.2, choose a $\nu$-measurable map $h(s)$ such that both $f_{t+h(s)}(s) = k(s) \widetilde{f}_{t+h(s)}(\Phi(s))$ and $f_{t+h(s)}(s) = a(s) f_t(s)$ a.e. $\lambda(dt)$ for $s \in \operatorname{proj}_S A$. Observe that $\operatorname{proj}_S A = C_P$ a.e. because $\operatorname{proj}_S \{(s,h): f_{t+h}(s) = a(h,s) f_t(s)$ a.e. $\lambda(dt)$ for some $a(h,s) \neq 0\} = C_P$ by the definition of $C_P$, and $\{(s,h): f_{t+h}(s) = k(s) \widetilde{f}_{t+h}(\Phi(s))$ a.e. $\lambda(dt)\} = S \times \mathbb{R}$ a.e. by (4.3). This then implies that, for a.e. $s \in C_P$, $a(s) k(s) \widetilde{f}_t(\Phi(s)) = a(s) f_t(s) = f_{t+h(s)}(s) = k(s) \widetilde{f}_{t+h(s)}(\Phi(s))$ a.e. $\lambda(dt)$ and hence that (4.5) holds.

Since (4.5) implies $s \in C_P \Rightarrow \Phi(s) \in \widetilde{C}_P$ a.e. $\nu(ds)$, we have $C_P \subset \Phi^{-1}(\widetilde{C}_P)$ a.e. $\nu(ds)$. Since $S = C_P$ a.e. $\nu(ds)$, we have

$$S = \Phi^{-1}(\widetilde{C}_P) \qquad \text{a.e. } \nu(ds).$$

This implies that $\widetilde{S} = \widetilde{C}_P$ a.e. $\widetilde{\nu}(d\widetilde{s})$. Indeed, if $\widetilde{\nu}(\widetilde{S} \setminus \widetilde{C}_P) > 0$, then by (4.4), we have $\nu(\Phi^{-1}(\widetilde{S} \setminus \widetilde{C}_P)) > 0$ as well. However, this contradicts $S = \Phi^{-1}(\widetilde{C}_P)$ a.e. $\nu(ds)$ since $\Phi^{-1}(\widetilde{C}_P)$ and $\Phi^{-1}(\widetilde{S} \setminus \widetilde{C}_P)$ are disjoint.

REMARK 4.1. In the case when $C_P$ is not equal to $S$ a.e. $\nu(ds)$, we may argue as above for the converse and show that if $\Phi(s) \in \widetilde{C}_P$, then $s \in C_P$ a.e. $\nu(ds)$. In other words, $\Phi^{-1}(\widetilde{C}_P) \subset C_P$ a.e. $\nu(ds)$. Since $C_P \subset \Phi^{-1}(\widetilde{C}_P)$ a.e. $\nu(ds)$ as shown above, we conclude that

$$(4.6) \qquad C_P = \Phi^{-1}(\widetilde{C}_P) \qquad \text{a.e. } \nu(ds).$$

Relationship (4.6) is used in the proof of Theorem 6.1.

CONTINUATION OF STEP 1. We may thus suppose without loss of generality that the representation $\{f_t\}_{t \in T}$ of $X_\alpha$ is minimal and that $C_P = S$ a.e. $\nu(ds)$. By Theorem 3.1 in Rosiński (1995), there is a flow $\{\phi_t\}_{t \in T}$ on $(S, \mathcal{S}, \nu)$ and a corresponding cocycle $\{a_t\}_{t \in T}$ such that, for all $t \in T$,

$$(4.7) \qquad f_t(s) = a_t(s) \left\{ \frac{d(\nu \circ \phi_t)}{d\nu}(s) \right\}^{1/\alpha} f_0(\phi_t(s))$$

a.e. $\nu(ds)$, where $f_0 \in L^\alpha(S, \mathcal{S}, \nu)$, that is, the process $X_\alpha$ is generated by the flow $\{\phi_t\}_{t \in T}$ in the sense of Definition 1.1.

STEP 2. We now show that the flow $\{\phi_t\}_{t \in T}$ is periodic. To do so, consider the set

$$(4.8) \quad A = \{(s,h) \in S \times (T \setminus \{0\}) : f_{t+h}(s) = a(h,s) f_t(s)$$
$$\text{a.e. } \lambda(dt) \text{ for some } a(h,s) \neq 0\}.$$

Observe now that by using (4.7) and the definition of a flow and a cocycle in Appendix A, for any $t, h \in T$,

$$f_{t+h}(s) = a_h(s) \left\{ \frac{d(\nu \circ \phi_h)}{d\nu}(s) \right\}^{1/\alpha} f_t(\phi_h(s))$$

a.e. $\nu(ds)$. Then, setting

$$(4.9) \quad A_0 = A \cap \{(s,h) \in S \times (T \setminus \{0\}) : f_{t+h}(s) = b(h,s) f_t(\phi_h(s))$$
$$\text{a.e. } \lambda(dt) \text{ for some } b(h,s) \neq 0\},$$

we have $A = A_0$ a.e. $\nu(ds)\tau(dh)$, where $\tau(dh)$ is any $\sigma$-finite measure on $T$. We now want to show that by setting

$$(4.10) \qquad A_1 = A_0 \cap \{(s,h) \in S \times (T \setminus \{0\}) : \phi_h(s) = s\},$$

we also have $A_1 = A_0$ a.e. $\nu(ds)\tau(dh)$. It is enough to prove that $\phi_h(s) = s$ a.e. $\nu(ds)\tau(dh)$ for $(s,h) \in A_0$. Supposing that this is not true, we can find $h$ such that $\phi_h(s) \neq s$ a.e. on a set of positive $\nu$ measure for $(s,h) \in A_0$ [otherwise, $\phi_h(s) = s$ a.e. $\nu(ds)\tau(dh)$ for $(s,h) \in A_0$ by the Fubini theorem]. Then, setting $\phi_h(s) = s$ for $(s,h) \in A_0^c$, we claim that a.e. $\lambda(dt)$,

$$(4.11) \qquad f_t(\phi_h(s)) = c(h,s) f_t(s) \qquad \text{a.e. } \nu(ds),$$

where $c(h,s) \neq 0$. This is clearly true for $(s,h) \in A_0^c$ since $\phi_h(s) = s$. This is also true for $(s,h) \in A_0$, because it follows from the definition of $A_0$ that the relationships $f_{t+h} = a f_t$ and $f_{t+h} = b f_t \circ \phi_h$ imply $f_t \circ \phi_h = c f_t$. We claim now that (4.11) is true not only a.e. $\lambda(dt)$, but for all $t \in T$. We only need to consider the case $T = \mathbb{R}$, because when $T = \mathbb{Z}$, the statements "a.e. $\lambda(dt)$" and "for all $t \in T$" are equivalent. Let $t \in \mathbb{R}$ be fixed. Since (4.11) holds a.e. $\lambda(dt)$, there is a sequence $\{t_n\}$ such that $t_n \to t$ and (4.11) holds with $t$ replaced by $t_n$. Since, by Lemma 4.3, $f_{t_n} \to f_t$ in $L^\alpha(S, \nu)$, we can select a subsequence $t_{n'}$ such that $f_{t_{n'}} \to f_t$ a.e. Then, relationship (4.11) with $t$ follows from analogous relationships with $t$ replaced by $t_{n'}$ by letting $n' \to \infty$. Together with the facts that $\phi_h(s) \neq s$ on a set of positive measure $\nu$ and $\phi_h(s)$ is nonsingular, (4.11) contradicts the minimality of the representation $\{f_t\}_{t \in T}$ [see Appendix B and, in particular, condition (M2$'$)]. Hence, $A_1 = A_0$ a.e. $\nu(ds)\tau(dh)$ and since $A_0 = A$ a.e. $\nu(ds)\tau(dh)$ as well, we have

$$(4.12) \qquad A = A_1 \qquad \text{a.e. } \nu(ds)\tau(dh).$$

By Lemma 4.2, we can choose a $\nu$-measurable function $h = h(s) \neq 0$ defined for $s \in \text{proj}_S A_1$ such that $(s, h(s)) \in A_1$ and, in particular,

$$(4.13) \qquad \phi_{h(s)}(s) = s.$$

By using (4.12), we have $\text{proj}_S A_1 = \text{proj}_S A$ a.e. $\nu(ds)$. Since $\text{proj}_S A = C_P$ by (4.1) and $C_P = S$ a.e. by assumption, we have $\text{proj}_S A_1 = S$ a.e. $\nu(ds)$, that is, (4.13) holds for a.e. $s \in S$. This shows that $S = P$ a.e. $\nu(ds)$, where $P$ are the periodic points of the flow $\{\phi_t\}_{t \in T}$ defined by (2.1).

STEP 3. We can now easily conclude the proof. We have

$$X_\alpha(t) \stackrel{d}{=} \int_P f_t(s) M_\alpha(ds) =: X_\alpha^P(t).$$

The process $X_\alpha^P$ is generated by a periodic flow $\{\phi_t|_P\}_{t \in T}$ and hence, by Theorem 3.1, it is a stationary periodic process. $\square$

EXAMPLE 4.1. Consider the process $X_\alpha$ of Example 3.2 defined through the kernel

$$f_t(v, z) = \cos(v + zt), \qquad v \in (0, 2\pi), z \in \mathbb{R} \text{ and } t \in \mathbb{R}.$$

Since $f_{t+2\pi z^{-1}}(v,z) = f_t(v,z)$ for all $v \in (0, 2\pi)$, $z \in \mathbb{R}$, the periodic component set $C_P$ associated with the representation $\{f_t\}_{t \in \mathbb{R}}$ is the full space $\mathbb{R} \times (0, 2\pi)$. Hence, Theorem 4.1 implies that, as already shown in Example 3.2, the process $X_\alpha$ is a stationary periodic process.

EXAMPLE 4.2. Consider the process

$$X_\alpha(t) = \int_0^1 \{v+t\}_1 M_\alpha(dv), \qquad t \in \mathbb{R}, \tag{4.14}$$

where $M_\alpha$ has control measure $dv$. Setting $f_t(v) = \{v+t\}_1$ for the kernel of $X_\alpha$, we see that $f_{t+1}(v) = f_t(v)$ for all $v \in [0, 1)$. Since the periodic component set $C_P$ associated with the representation $\{f_t\}$ is the full space $[0, 1)$, $X_\alpha$ is a stationary periodic process. This fact can also be seen directly from the representation (3.2) which becomes that of $X_\alpha$ when $Z = \{1\}$, $b_1(1) = 1$, $q(1) = 1$ and $g(1, v) = v$.

Observe also that the representation (4.14) is minimal (see Appendix B). Indeed, taking $t \in (0, 1)$, since $f_0(v) < f_{1-t}(v)$ for $v \in [0, t)$ and $f_{1-t}(v) < f_0(v)$ for $v \in [t, 1)$, we have that $(f_{1-t}/f_0)^{-1}([1, \infty)) = [0, t)$ and hence that $[0, t) \in \sigma\{f_t/f_s, s, t \in \mathbb{R}\}$. It follows that $\mathcal{B}([0, 1)) = \sigma\{f_t/f_s, s, t \in \mathbb{R}\}$. The condition $\mathrm{supp}\{f_t, t \in \mathbb{R}\} = [0, 1)$ is obviously satisfied.

Since the kernel of the stationary periodic process $X_\alpha$ in Example 4.2 satisfies relationship $f_{t+1}(v) = f_t(v)$ for all $v \in [0, 1), t \in \mathbb{R}$, we have $X_\alpha(t) = X_\alpha(t+1)$ for all $t \in \mathbb{R}$, and hence the process $X_\alpha$ is not ergodic. Since the kernels of other stationary periodic processes are also periodic in nature, we can expect that all these processes are not ergodic either. The following theorem shows that this is indeed the case.

THEOREM 4.2. *The $S\alpha S$, $\alpha \in (0, 2)$, stationary periodic processes are not ergodic.*

PROOF. In view of Corollary 3.1, we can suppose that a stationary periodic process has a representation (3.2). We may suppose without loss of generality that the measure $\sigma(dz)\lambda(dv)$ is finite on $S := Z \times [0, q(\cdot))$, because, otherwise, we may replace $\sigma(dz)\lambda(dv)$ by a finite measure $k(z)^\alpha \sigma(dz)\lambda(dv)$, where $k(z) > 0$ satisfies

$$\int_Z \int_{[0,q(z))} k(z)^\alpha \sigma(dz)\lambda(dv) = \int_Z k(z)^\alpha \lambda([0, q(z)))\sigma(dz) < \infty,$$

and define the process $X_\alpha$ as in (3.2) with $g(z, \{v+t\}_{q(z)})$ divided by $k(z)$. Then a stationary periodic process (3.2) is generated by a flow $\phi_t(z, v) = (z, \{v+t\}_{q(z)})$ on a space $S = Z \times [0, q(\cdot))$ such that, without loss of generality, $(\sigma \otimes \lambda)(S) < \infty$. Observe that the measure $\nu := \sigma \otimes \lambda$ is invariant under the flow $\phi_t$ since $d(\nu \circ \phi_t)/d\nu = 1$. Since $\nu(S \cap \phi_t S) = \nu(S)$ is not only finite but also positive and does not depend on $t \in T$, we have $\lim_{t \to \infty} \nu(S \cap \phi_t S) \neq 0$. Applying Theorem 4.1 in Gross (1994) [see also Corollary 2.1 of Rosiński and Samorodnitsky (1996)], we conclude that a stationary periodic process is not mixing. Applying the same result of Gross together with a statement at the top of page 279 in Gross (1994) [see also Remark 2.3 of Rosiński and Samorodnitsky (1996)], it follows that it is not weak-mixing either and since, for stable processes, weak mixing and ergodicity coincide [see Samorodnitsky and Taqqu (1994), page 580], it follows that the process is not ergodic. □

Finally, we establish the results used in the proofs of Lemma 4.1 and Theorem 4.1.

LEMMA 4.2. *Let $(S_1, \mathcal{S}_1, \nu_1)$ and $(S_2, \mathcal{S}_2, \nu_2)$ be two standard Lebesgue spaces and let $(S_1 \times S_2, \mathcal{S}_1 \otimes \mathcal{S}_2, \nu_1 \otimes \nu_2)$ be their Cartesian product. Let also $A \in \mathcal{S}_1 \otimes \mathcal{S}_2$ be a Borel set of $S_1 \times S_2$. Then the set*

$$\operatorname{proj}_{S_1} A := \{s_1 \in S_1 : \exists\, s_2 \in S_2 : (s_1, s_2) \in A\}$$

*is $\nu_1$-measurable and there is a $\nu_1$-measurable function $h : \operatorname{proj}_{S_1} A \mapsto A$ such that $(s_1, h(s_1)) \in A$ for all $s_1 \in \operatorname{proj}_{S_1} A$.*

PROOF. The set $\operatorname{proj}_{S_1} A$ is $\nu_1$-measurable because the map $\operatorname{proj}_{S_1}(s_1, s_2) = s_1$ is continuous and the set $A$ can be approximated $(\nu_1 \otimes \nu_2)$-a.e. by rectangles whose projections are measurable. We show next that there is a $\nu_1$-measurable map $h : \operatorname{proj}_{S_1} A \mapsto A$ such that $(s_1, h(s_1)) \in A$ for $s_1 \in \operatorname{proj}_{S_1} A$. To do so, we use Theorem 3.4.3 in Arveson (1976), page 77, which concerns the so-called cross sections of Borel maps. Consider the map $f = \operatorname{proj}_{S_1} : A \mapsto f(A) = \operatorname{proj}_{S_1} A$. The image set $f(A)$, together with the induced Borel structure $\mathcal{F}(A) = \{f(A) \cap B : B \in \mathcal{S}_1\}$, is a Borel space. Moreover, this Borel space is countably separated [as defined in Arveson (1976), page 69] since the underlying standard Lebesgue space $(S_1, \mathcal{S}_1, \nu_1)$ is countably separated. The Borel set $A$, equipped with the Borel structure $\mathcal{A} = \{A \cap B : B \in \mathcal{S}_1 \otimes \mathcal{S}_2\}$, is also a Borel space. It is an analytic Borel space [as defined in Arveson (1976), page 71] by using the Corollary in Arveson (1976), page 65, and the fact that $A$ is a Borel set. Since $f^{-1}(f(A) \cap B) = A \cap (B \times \mathbb{R}) \in \mathcal{A}$ for all $B \in \mathcal{S}_1$, the map $f : (A, \mathcal{A}) \mapsto (f(A), \mathcal{F}(A))$ is Borel. It follows from Theorem 3.4.3 in Arveson (1976) that there is a $\nu_1$-measurable map $g : f(A) \mapsto A$ such that $f(g(s_1)) = s_1$. Since $f$ is a projection, we have that $g(s_1) = (s_1, h(s_1))$ for some $\nu$-measurable map $h(s_1)$ and hence that there is a $\nu$-measurable map $h(s_1)$ such that $(s_1, h(s_1)) \in A$. □

The next result follows from Surgailis, Rosiński, Mandrekar and Cambanis (1998), who considered measurable stationary increments processes. We present their proof here for the convenience of the reader.

LEMMA 4.3. *Every measurable stationary process $\{X(t)\}_{t\in\mathbb{R}}$ is continuous in probability.*

PROOF. Consider $\|\xi\| = E\min(|\xi|, 1)$ defined for random variables $\xi \in L^0(\Omega, P)$. (We use $\|\cdot\|$ as a convenient notation.) To show that a process $\{X(t)\}_{t\in\mathbb{R}}$ is continuous in probability, it is enough to prove that $\|X(t) - X(s)\| \to 0$ as $s \to t$.

Fix an arbitrary small $\varepsilon > 0$. For $t \in \mathbb{R}$, let $B_t = \{s \in \mathbb{R} : \|X(t) - X(s)\| < \varepsilon\}$. Since the process $X$ is assumed measurable, by Theorem 3 in Cohn (1972), the map $\mathbb{R} \ni t \mapsto X(t)$ is Borel and has a separable range. Hence, we can choose a sequence $\{t_n\} \subset \mathbb{R}$ such that $\{B_{t_n}\}$ are Borel and $\mathbb{R} = \bigcup_n B_{t_n}$. Then, there is $t_n$ such that the Lebesgue measure of $B_{t_n}$ is positive. By the Steinhaus lemma [see, e.g., Bingham, Goldie and Teugels (1987), page 2], the set $B = B_{t_n} - B_{t_n}$ (of points $z$ such that $z = x - y$, $x, y \in B_n$) contains an open interval $(-\delta, \delta)$ with some $\delta > 0$. If $|s - t| < \delta$, then $s - t = u - v$ for some $u, v \in B_{t_n}$, and hence

$$\|X(t) - X(s)\| = \|X(u) - X(v)\| \leq \|X(u) - X(t_n)\| + \|X(v) - X(t_n)\| < 2\varepsilon,$$

where we used the stationarity of $X$ and the fact that $u, v \in B_{t_n}$. This shows that $\|X(t) - X(s)\| \to 0$ as $s \to t$. □

**5. Characterization of stationary cyclic processes.** We know from Section 4 how to identify stationary periodic processes. We want to identify stationary cyclic processes, namely stationary periodic processes without a harmonizable (or trivial) component (see Definition 3.2). Rosiński (1995) showed that harmonizable processes (or trivial processes in the real-valued case) can be identified through the *harmonizable (or trivial) component set*

$$(5.1) \quad C_F = \{s \in S : f_{t_1+t_2}(s)f_0(s) = f_{t_1}(s)f_{t_2}(s) \text{ a.e. } \lambda(dt_1)\lambda(dt_2)\}.$$

LEMMA 5.1. *We have*

$$(5.2) \quad C_F \subset C_P \quad a.e. \ \nu(ds).$$

PROOF. By Lemma 5.5 in Rosiński (1995), $f_0 \neq 0$ a.e. on $C_F$. Hence, by fixing $t_1 = h$ in the definition (5.1) of $C_F$, we get that, for a.e. $s \in C_F$, $f_{t+h}(s) = a(s)f_t(s)$ a.e. $\lambda(dt)$ with $a(s) = f_h(s)/f_0(s)$. This shows (5.2). □

Since stationary cyclic processes are stationary periodic processes without a harmonizable (or trivial) component, we expect that stationary cyclic processes can be identified through the set $C_L = C_P \setminus C_F$.

DEFINITION 5.1. *A cyclic component set* of a stationary stable process $X_\alpha$ having a representation (1.1) is defined as

$$C_L = C_P \setminus C_F, \tag{5.3}$$

where the sets $C_P$ and $C_F$ are defined by (4.1) and (5.1), respectively.

The next result shows that stationary cyclic processes can indeed be identified through the cyclic component set $C_L$.

THEOREM 5.1. *A $S\alpha S$, $\alpha \in (0,2)$, stationary process $X_\alpha$ given by* (1.1) *with* $\mathrm{supp}\{f_t, t \in T\} = S$ *a.e.* $\nu(ds)$ *is a stationary cyclic process if and only if*

$$C_L = S \qquad a.e.\ \nu(ds),$$

*where $C_L$ is the cyclic component set defined in* (5.3).

PROOF. If $X_\alpha$ is a stationary cyclic process, then it is a stationary periodic process as well. However, by Theorem 4.1, $C_P = S$ a.e. $\nu(ds)$. By (5.3), we have $C_P = C_F + C_L$. Since $X_\alpha$ does not have a harmonizable (or trivial) component, Rosiński (1995) results show that $C_F = \varnothing$ a.e. $\nu(ds)$ and hence that $C_L = S$ a.e. $\nu(ds)$. Conversely, if $C_L = S$ a.e. $\nu(ds)$, then $C_F = \varnothing$ a.e. $\nu(ds)$ and $C_P = S$ a.e. $\nu(ds)$. Hence, by Theorem 4.1, the process $X_\alpha$ is a stationary periodic process. Since $C_F = \varnothing$, by Rosiński (1995), the process $X_\alpha$ does not have a harmonizable (or trivial) component. □

Since $C_F \subset C_P$ by Lemma 5.1, we may ask how the definition (5.3) of the set $C_F$ relates to the definition (4.1) of the periodic component set $C_P$. The following result provides an answer.

PROPOSITION 5.1. *We have*

$$\begin{aligned}C_F = \{s \in S : \exists T \setminus \{0\} \ni h_n = h_n(s) \to 0 \text{ as } n \to \infty : \\ f_{t+h_n}(s) = a(h_n, s) f_t(s) \text{ a.e. } \lambda(dt) \text{ for some } a(h_n, s) \neq 0\}\end{aligned} \tag{5.4}$$

*a.e.* $\nu(ds)$ *when* $T = \mathbb{R}$ *and*

$$C_F = \{s \in S : f_{t+1}(s) = a(s) f_t(s) \text{ a.e. } \lambda(dt) \text{ for some } a(s) \neq 0\} \tag{5.5}$$

*a.e.* $\nu(ds)$ *when* $T = \mathbb{Z}$.

PROOF. We first consider the case $T = \mathbb{R}$. Denote the set on the right-hand side of (5.4) by $C_F^0$. Let us first show that $C_F \subset C_F^0$ a.e. $\nu(ds)$. As shown in the proof of Theorem 5.7 in Rosiński (1995), for each $t \in T$, $f_t(s) = e^{itk(s)} f_0(s)$ a.e. $s \in C_F$, where $k(s)$ is some function [in the real-valued case, the relationship is $f_t(s) = f_0(s)$]. By the Fubini's theorem, we also have that, for a.e. $s \in C_F$, $f_t(s) = e^{itk(s)}f_0(s)$ a.e. $\lambda(dt)$. Then, since for a.e. $s \in C_F$, $f_{t+h}(s) = e^{ihk(s)}f_t(s)$ a.e. $\lambda(dt)$ for any $h \in \mathbb{R}$, it holds in particular for a sequence $h_n$ (not even depending on $s$) satisfying $h_n \to 0$. Setting $a(h_n, s) = e^{ih_n k(s)}$, we obtain $C_F \subset C_F^0$ a.e. $\nu(ds)$.

We now show that $C_F^0 \subset C_F$ a.e. $\nu(ds)$. Let $\{\widetilde{f}_t\}_{t \in T}$ be the kernel in a minimal integral representation (4.2) for the process $X_\alpha$. Let also $\widetilde{C}_F$ and $\widetilde{C}_F^0$ be the sets defined in the same way as $C_F$ and $C_F^0$, but by using only the kernel $\widetilde{f}_t$. We can show as in the proof of Theorem 4.1 that $C_F^0 = \Phi^{-1}(\widetilde{C}_F^0)$ a.e. $\nu(ds)$, where $\Phi$ is the map appearing in (4.3). Moreover, as shown in the proof of Theorem 5.7 in Rosiński (1995), $C_F = \Phi^{-1}(\widetilde{C}_F)$ a.e. $\nu(ds)$. It is then enough to show that $\widetilde{C}_F^0 \subset \widetilde{C}_F$ a.e. or, equivalently, that $C_F^0 \subset C_F$ a.e., but where $\{f_t\}_{t \in T}$ is supposed to be a minimal representation. If $\{f_t\}_{t \in T}$ is minimal, then it is generated by a flow $\{\phi_t\}_{t \in T}$ in the sense of Definition 1.1 [Theorem 3.1 in Rosiński (1995)]. By Lemma 5.2, the set $C_F^0$ is a.e. invariant under the flow $\{\phi_t\}_{t \in T}$. The process

$$\int_{C_F^0} f_t(s) M_\alpha(ds)$$

is then stationary, its representation $\{f_t|_{C_F^0}\}_{t \in T}$ is minimal and is generated by the flow $\{\phi_t|_{C_F^0}\}_{t \in T}$. Arguing as in the proof of Theorem 4.1 [see (4.13)], we can show that, for a.e. $s \in C_F^0$,

(5.6) $$\phi_{h_n(s)}(s) = s \quad \text{for } T \setminus \{0\} \ni h_n(s) \to 0.$$

In view of the special representation (A.2) of a flow without fixed points, the last relationship cannot hold for points which are not fixed and hence we obtain that, for a.e. $s \in C_F^0$, $\phi_t(s) = s$ for all $t \in T$. Then, by Proposition 5.8 in Rosiński (1995), $C_F^0 \subset C_F$ a.e. $\nu(ds)$.

The case $T = \mathbb{Z}$ can be proved in a similar way. The main difference is that (5.6) is replaced by $\phi_1(s) = s$ for a.e. $s \in C_F^0$. This shows that $\phi_t(s) = s$ a.e. $s \in C_F$ for all $t \in T$ and hence $C_F^0 \subset C_F$ a.e. $\nu(ds)$ as well. $\square$

The sets $C_F$ and $C_P$ are explicitly identified by (5.1) and (4.1), respectively. Proposition 5.1 yields the following explicit identification of $C_L = C_P \setminus C_F$.

COROLLARY 5.1. *We have*

(5.7) $$C_L = \{s \in S : \exists h_0 = h_0(s) \in T \setminus \{0\}, \not\exists T \setminus \{0\} \ni h_n = h_n(s) \to 0 \text{ as } n \to \infty :$$
$$f_{t+h_n}(s) = a(h_n, s) f_t(s) \text{ a.e. } \lambda(dt), n \geq 0, \text{ for some } a(h_n, s) \neq 0\}$$

a.e. $\nu(ds)$ when $T = \mathbb{R}$ and

$$C_L = \{s \in S : \exists h = h(s) \in T \setminus \{0\} : f_{t+h}(s) = a(h,s)f_t(s)$$

(5.8) $$\text{a.e. } \lambda(dt) \text{ for some } a(h,s) \neq 0\}$$

$$\cap \{s \in S : f_{t+1}(s) \neq a(s)f_t(s) \text{ a.e. } \lambda(dt) \text{ for all } a(s) \neq 0\}$$

a.e. $\nu(ds)$ when $T = \mathbb{Z}$.

EXAMPLE 5.1. The real part of a harmonizable process $X_\alpha$ of Examples 3.2 and 4.1 is a stationary cyclic process, because $C_P = \mathbb{R} \times [0, 2\pi)$ as shown in Example 4.1 and $C_F = \varnothing$ a.e. by using Proposition 5.1. To see that $C_F = \varnothing$ a.e., observe that the condition $f_{t+h}(s) = a(h,s)f_t(s)$ a.e. $dt$ for the process $X_\alpha$ becomes $\cos(v + z(t+h)) = a(h,z,v)\cos(v + zt)$ a.e. $dt$. After fixing $v$ and $z \neq 0$, we get $\cos(w + zh) = a(h,z)\cos(w)$ a.e. $dw$. This holds only for $h = \pi k/z \neq 0$, $k \in \mathbb{Z} \setminus \{0\}$ [with $a(h,z) = (-1)^k$], which cannot be taken arbitrarily small.

The process of Example 4.2 is also a stationary cyclic process since $C_F = \varnothing$ as implied by Proposition 5.1.

Finally, we establish an auxiliary result used in the proof of Proposition 5.1.

LEMMA 5.2. *If $\{f_t\}_{t \in T}$ is a representation of a $S\alpha S$, $\alpha \in (0,2)$, stationary process generated by a flow $\{\phi_t\}_{t \in T}$ in the sense of Definition 1.1 and $C_F^0$ denotes the set on the right-hand side of either (5.4) or (5.5), then $C_F^0$ is a.e. invariant under the flow $\{\phi_t\}_{t \in T}$, that is, $\nu(C_F^0 \triangle \phi_t^{-1}(C_F^0)) = 0$ for all $t \in T$.*

PROOF. We have to show that for $t_0 \in T$, $C_F^0 = \phi_{t_0}^{-1}(C_F^0)$ a.e., but since the flow $\{\phi_t\}_{t \in T}$ satisfies the group property, it is enough to show that for $t_0 \in T$, $C_F^0 \subset \phi_{t_0}^{-1}(C_F^0)$, that is, $s \in C_F^0$ implies $\phi_{t_0}(s) \in C_F^0$ a.e. $\nu(ds)$. We consider only the case $T = \mathbb{R}$. Observe that by using (1.4), for $t, h \in T$,

(5.9) $$f_{t+t_0+h}(s) = a_{t_0}(s) \left\{ \frac{d(\nu \circ \phi_{t_0})}{d\nu}(s) \right\}^{1/\alpha} f_{t+h}(\phi_{t_0}(s))$$

a.e. $\nu(ds)$. By Lemma 4.2 above, we can choose a sequence of $\nu$-measurable functions $h_n : C_F^0 \mapsto T \setminus \{0\}$ such that $h_n(s) \to 0$ and, for $s \in C_F^0$, $f_{t+t_0+h_n(s)}(s) = a_n(s)f_{t+t_0}(s)$ a.e. $\lambda(dt)$ for some $a_n(s) \neq 0$, and that for a.e. $s \in C_F^0$, the relationship (5.9) holds a.e. $\lambda(dt)$ with $h$ replaced by $h_n(s)$. Since we also have that by (5.9), for a.e. $\nu(ds)$,

$$f_{t+t_0}(s) = a_{t_0}(s) \left\{ \frac{d(\nu \circ \phi_{t_0})}{d\nu}(s) \right\}^{1/\alpha} f_t(\phi_{t_0}(s))$$

a.e. $\lambda(dt)$, it follows that, for a.e. $s \in C_F^0$, $f_{t+h_n(s)}(\phi_{t_0}(s)) = b_n(s)f_t(\phi_{t_0}(s))$ a.e. $\lambda(dt)$ for some $b_n(s) \neq 0$, that is, $\phi_{t_0}(s) \in C_F^0$ a.e. $\nu(ds)$. □

**6. Further decomposition of stationary stable processes.** In this section, we refine the decomposition (1.13) of Rosiński (1995) by showing that (1.18) holds. We first need to recall some basic facts behind the decompositions (1.6) and (1.13) found in Rosiński (1995). Consider a S$\alpha$S stationary process $X_\alpha$ given by (1.1) with $\mathrm{supp}\{f_t, t \in T\} = S$ a.e. $\nu(ds)$. Let

$$\text{(6.1)} \qquad D = \left\{ s \in S : \int_T |f_t(s)|^\alpha \lambda(dt) < \infty \right\}$$

and

$$\text{(6.2)} \qquad C = \left\{ s \in S : \int_T |f_t(s)|^\alpha \lambda(dt) = \infty \right\}.$$

If, in addition, the process $X_\alpha$ is generated by a flow $\{\phi_t\}_{t \in T}$ in the sense of Definition 1.1, then $D$ and $C$ are (a.e.) the dissipative and the conservative parts of the flow $\{\phi_t\}_{t \in T}$, respectively [see Rosiński (1995)]. Recall also the definition (5.1) of a harmonizable (or trivial) component set $C_F$. If the process $X_\alpha$ is generated by a flow $\{\phi_t\}_{t \in T}$ and if its representation is minimal (in the sense of Appendix B), then $C_F = F$ a.e., where $F$ is the set of the fixed points defined by (2.2). We can show that $C_F \subset C$ a.e. [see Rosiński (1995)].

The processes $X_\alpha^D$ and $X_\alpha^C$ in (1.6) are then defined as

$$\text{(6.3)} \qquad X_\alpha^D(t) = \int_D f_t(s) M_\alpha(ds), \qquad X_\alpha^C(t) = \int_C f_t(s) M_\alpha(ds),$$

that is, by replacing the full space $S$ in (1.1) by its disjoint subsets $D$ and $C$, respectively. The other two processes on the right-hand side of (1.13) are defined as

$$\text{(6.4)} \qquad X_\alpha^F(t) = \int_{C_F} f_t(s) M_\alpha(ds), \qquad X_\alpha^{C \setminus F}(t) = \int_{C \setminus C_F} f_t(s) M_\alpha(ds).$$

As one can see from these definitions, the idea behind decompositions (1.6) and (1.13) is to partition the underlying space $S$ into appropriately chosen subsets and then define the processes in a decomposition as integrals over these subsets. To get the decomposition (1.18), we pursue the same idea. We use the following lemma.

LEMMA 6.1. *We have*

$$\text{(6.5)} \qquad C_P \subset C \qquad a.e. \ \nu(ds),$$

*where the sets $C_P$ and $C$ are defined by* (4.1) *and* (6.2).

PROOF. Observe that, for $s \in C_P$,

$$\int_T |f_t(s)|^\alpha \lambda(dt) = \int_{[0,|h(s)|)} |f_t(s)|^\alpha \lambda(dt) \left( \sum_{n=-\infty}^{\infty} |a(h,s)|^{\alpha n} \right) = \infty.$$

Hence $s \in C$ by (6.2). □

By using Lemmas 5.1 and 6.1, and the definition of $X_\alpha^{C \setminus F}$ in (6.4), we can write

$$X_\alpha^{C \setminus F} \stackrel{d}{=} X_\alpha^L + X_\alpha^{C \setminus P}, \tag{6.6}$$

where

$$(6.7) \quad X_\alpha^L(t) = \int_{C_L} f_t(s) M_\alpha(ds), \qquad X_\alpha^{C \setminus P}(t) = \int_{C \setminus C_P} f_t(s) M_\alpha(ds).$$

The following result refines the decomposition (1.13) of Rosiński (1995). We say that the decomposition (6.6) is unique in distribution if the distribution of its two components does not depend on the representation (1.1) of the process $X_\alpha$. We also say that a stable stationary process does not have a periodic component if it cannot be represented as the sum of two independent S$\alpha$S processes, one of which is a nondegenerate stationary periodic process.

THEOREM 6.1. *The decomposition* (6.6) *is unique in distribution. Moreover, the process* $X_\alpha^L$ *is a stationary cyclic process. The process* $X_\alpha^{C \setminus P}$ *is a S$\alpha$S stationary process generated by a conservative flow without a periodic component.*

PROOF. The idea is similar to that of the proof of Theorem 5.7 in Rosiński (1995). Let $\{\widetilde{f}_t\}_{t \in T}$ be a minimal integral representation (4.2) for the process $X_\alpha$ so that, in particular, relationships (4.3) and (4.4) hold. Let also $\widetilde{C}$, $\widetilde{C}_F$, $\widetilde{C}_L$ and $\widetilde{C}_P$ be the sets defined in (6.2), (5.1), (5.3) and (4.1), respectively, by using the functions $\widetilde{f}_t$. Let also $\widetilde{X}_\alpha^L$ and $\widetilde{X}_\alpha^{C \setminus P}$ be the two components on the right-hand side of the decomposition (6.6) for the process $X_\alpha$ obtained by using the sets $\widetilde{C}, \widetilde{C}_F, \widetilde{C}_L, \widetilde{C}_P$ and the kernel $\widetilde{f}_t$. It is enough to show that $X_\alpha^L \stackrel{d}{=} \widetilde{X}_\alpha^L$ and $X_\alpha^{C \setminus P} \stackrel{d}{=} \widetilde{X}_\alpha^{C \setminus P}$. As shown in the proofs of Theorems 4.3 and 4.7 in Rosiński (1995) and in (4.6) in Theorem 4.1, we have $C = \Phi^{-1}(\widetilde{C})$, $C_F = \Phi^{-1}(\widetilde{C}_F)$ and $C_P = \Phi^{-1}(\widetilde{C}_P)$ a.e. $\nu(ds)$, where $\Phi$ is the map appearing in (4.3). Hence, $C_L = C_P \setminus C_F = \Phi^{-1}(\widetilde{C}_P \setminus \widetilde{C}_F) = \Phi^{-1}(\widetilde{C}_L)$ and $C \setminus C_P = \Phi^{-1}(\widetilde{C} \setminus \widetilde{C}_P)$ a.e. $\nu(ds)$. The two required identities in distribution above follow from these relationships by using (4.3), (4.4) and a change of variables as at the end of the proof of Theorem 4.3 in Rosiński (1995).

Let us show now that the process $X_\alpha^L$ is a stationary cyclic process. By using the preceding discussion, we may suppose without loss of generality that the representation $\{f_t\}_{t \in T}$ for the process $X_\alpha$ is minimal and hence, by Theorem 3.1 in Rosiński (1995), generated by a flow $\{\phi_t\}_{t \in T}$ in the sense of Definition 1.1. By using Lemma 6.2, the set $C_P$ is a.e. invariant under the flow $\{\phi_t\}_{t \in T}$. Since the set $C_F$ is a.e. invariant under the flow by Lemma 5.6 in Rosiński (1995), the set $C_L = C_P \setminus C_F$ is a.e. invariant under the flow as well. Consequently, the process $X_\alpha^L$ defined on the set $C_L$ is stationary. It is a stationary cyclic process by construction in view of Theorem 5.1.

We now focus on the process $X_\alpha^{C \setminus P}$. Since $X_\alpha^L$ and $X_\alpha^{C \setminus P}$ are independent and since $X_\alpha^{C \setminus F}$ and $X_\alpha^L$ are both S$\alpha$S, stationary and conservative, so is the process $X_\alpha^{C \setminus P}$ ["conservative" follows from the uniqueness in distribution of the decomposition (1.6)]. It remains to show that $X_\alpha^{C \setminus P}$ does not have a periodic component. We do so by adapting the end of the proof of Theorem 5.7 in Rosiński (1995) to our case. Suppose that $X_\alpha^{C \setminus P}$ admits a periodic component, that is,

$$X_\alpha^{C \setminus P} \stackrel{d}{=} V + W,$$

where $V$ and $W$ are independent S$\alpha$S processes and $W$ is a nondegenerate stationary periodic process. Let $\{f_t^{C \setminus P}\}_{t \in T}$ be the representation of the process $X_\alpha^{C \setminus P}$, that is, the restriction of $f_t$ to $C \setminus C_P$, and let $\{g_t\}_{t \in T}$ be the representation of the process $W$ defined in (3.2) on the space $(Z \times [0, q(\cdot)), \sigma(dz)\lambda(dv))$ by using functions $b_1(z), q(z)$ and $g(z, v)$. By using Theorem 1.1 in Rosiński (1995), we obtain that

(6.8) $\qquad g_t(z, v) = h(z, v) f_t^{C \setminus P}(\Psi(z, v)) \qquad$ a.e. $\lambda(dt)\sigma(dz)\lambda(dv),$

where $\Psi : Z \times [0, q(\cdot)) \mapsto C \setminus C_P$ and $h : Z \times [0, q(\cdot)) \mapsto \mathbb{R}$ or $\mathbb{C}$ are some maps. Since $\sigma$ is not a zero measure (otherwise $W$ would be degenerate), then (6.8) is a contradiction to the fact that $\Psi(z, u) \in C \setminus C_P$ in view of the definitions of a stationary periodic process and $C_P$. $\square$

The following result was used in the proof of Theorem 6.1.

LEMMA 6.2. *If $\{f_t\}_{t \in T}$ is a representation of a S$\alpha$S, $\alpha \in (0, 2)$, stationary process generated by a flow $\{\phi_t\}_{t \in T}$ in the sense of Definition 1.1, then the periodic component set $C_P$ in (4.1) is a.e. invariant under the flow $\{\phi_t\}_{t \in T}$, that is, $\nu(C_P \triangle \phi_t^{-1}(C_P)) = 0$ for all $t \in T$.*

PROOF. The proof of this lemma is similar and simpler than that of Lemma 5.2, and hence is omitted. $\square$

The real part of a harmonizable process in (3.13) and the process in Example 4.2 are examples of stationary cyclic processes $X_\alpha^L$ in the decomposition (6.6) (see Example 5.1). We now provide examples of the "fourth" kind of processes $X_\alpha^{C \setminus P}$ in that decomposition. We consider the case $T = \mathbb{R}$ only. Extensions to $T = \mathbb{Z}$ are elementary.

EXAMPLE 6.1. Let $\{Y(t)\}_{t\in\mathbb{R}}$ be a stationary process which has càdlàg (that is, right-continuous and with limits from the left) paths, satisfies $P(|Y(t)| < c) < 1$ for all $c > 0$, $E|Y(t)|^\alpha < \infty$ and is ergodic. Let also $\Omega = \{w\}$ be the space of càdlàg functions on $\mathbb{R}$ in the real-valued case, let $\Omega = \{w = w_1 + iw_2\}$ be the space of càdlàg functions $w_1$ and $w_2$ on $\mathbb{R}$ in the complex-valued case, and let $P(dw)$ be a probability measure on the space $\Omega$ corresponding to the process $Y$. Consider now the process

$$(6.9) \qquad X_\alpha(t) = \int_\Omega f_t(w) M_\alpha(dw), \qquad t \in \mathbb{R},$$

where

$$f_t(w) = w(t) \quad \text{and} \quad M_\alpha(dw) \text{ has the control measure } P(dw).$$

We can show that $X_\alpha$ is a well-defined S$\alpha$S stationary process. When $Y$ is a real-valued Gaussian process, the process (6.9) is called a S$\alpha$S sub-Gaussian stationary process [see, e.g., Samorodnitsky and Taqqu (1994)]. Ergodicity of $Y$ is equivalent to the continuity of its spectral measure [see page 163 in Rozanov (1967)]. When $Y$ is a S$\alpha'$S stationary process with $\alpha < \alpha'$, the process (6.9) is called a substable process [page 143 in Samorodnitsky and Taqqu (1994)]. To prove that $X_\alpha$ is indeed the "fourth kind" process in the decomposition (6.6), it is enough to show that $X_\alpha$ is generated by a conservative flow and that $C_P = \varnothing$ a.e. $dP$.

Since $\Omega$, equipped with the usual Skorokhod $J_1$-topology, is a complete separable metric space, the space $(\Omega, \mathcal{F}, P)$, where $\mathcal{F}$ is the $\sigma$-field of the Borel sets, is standard Lebesgue. Observe now that

$$f_t(w) = w(t) = (\phi_t(w))(0) = f_0(\phi_t(w)),$$

where $\phi_t : \Omega \mapsto \Omega$ is defined by $(\phi_t(w))(s) = w(t+s)$. The collection of maps $\{\phi_t\}_{t\in\mathbb{R}}$ is a measurable flow on a standard Lebesgue space $\Omega$ and hence, in view of Definition 1.1, the process $X_\alpha$ is generated by the flow $\{\phi_t\}_{t\in\mathbb{R}}$. The flow is conservative since it is measure preserving and the measure $P$ on $\Omega$ is finite (in other words, there can be no wandering set of positive measure).

Let us show now that $C_P = \varnothing$ a.e. $dP$. By the definition (4.1) of $C_P$, we have

$$C_P = \{w \in \Omega : \exists\, h = h(w) \neq 0 : w(t+h) = a(h,w)w(t)$$

$$\text{a.e. } dt \text{ for some } a(h,w) \neq 0\}.$$

If $w \in C_P$ and $|a| \neq 1$, then $w(t) \to 0$ when either $t \to +\infty$ a.e. $dt$ or $t \to -\infty$ a.e. $dt$ ($t \to +\infty$ a.e. $dt$, for example, means that $t \to +\infty$ on a set $B$ such that $B^c = \varnothing$ a.e. $dt$). In either case, the $P$ measure of such sets is zero. For example, if $w(t) \to 0$ as $t \to +\infty$ a.e. $dt$, then $T^{-1} \int_0^T \mathbb{1}_{\{|w(t)|<1\}}(t)\, dt \to 1$ as $T \to \infty$. If the $P$ measure of that set is positive, this would contradict the

STABLE STATIONARY PROCESSES 35ergodicity according to which the limit is $P(|w(0)| < 1) < 1$. If $w \in C_P$ and $|a| = 1$, then $w$ is bounded a.e. $dt$ on $\mathbb{R}$. Supposing that $P(w \in C_P, |a| = 1) > 0$, we obtain a contradiction in the same way as above by considering the integral $T^{-1} \int_0^T \mathbb{1}_{\{|w(t)| < N\}}(t)\, dt$ for large enough $N$. Hence, $P(C_P) = 0$.

If we work exclusively with minimal representations of stationary stable processes, then we can relate the set $C_L = C_P \setminus C_F$ used in the definition (6.7) of the process $X_\alpha^L$ to the set of cyclic points $L$ of the underlying flow. This extends to the cyclic case, Proposition 5.8, in Rosiński (1995), where the set $C_F$ is identified as the set of the fixed points $F$ of the flow.

PROPOSITION 6.1. *If the representation $\{f_t\}_{t \in T}$ is minimal for the process $X_\alpha$ and $\{\phi_t\}_{t \in T}$ is the flow related to $\{f_t\}_{t \in T}$ in the sense of Definition* 1.1, *then*

(6.10) $$C_L = L \qquad a.e.\ \nu(ds),$$

*where $L$ is the set of cyclic points of the flow $\{\phi_t\}_{t \in T}$ defined in* (2.3).

PROOF. By Proposition 5.8 in Rosiński (1995), $C_F = F$ a.e., where $F$ is the set of the fixed points of the flow. It is therefore enough to prove that $C_P = P$ a.e. By using the fact that $\{f_t|_{C_P}\}_{t \in T}$ is a minimal representation for a S$\alpha$S process generated by the flow $\{\phi_t|_{C_P}\}_{t \in \mathbb{R}}$ and arguing as in the proof of Theorem 4.1, we can deduce [see (4.13)] that $\phi_{h(s)}(s) = s$ for a.e. $s \in C_P$. This shows that $C_P \subset P$ a.e. Suppose, on the other hand, that $\nu(P \setminus C_P) > 0$ and consider the process

(6.11) $$\int_{P \setminus C_P} f_t(s) M_\alpha(ds).$$

Since $P$ and $C_P$ are a.e. invariant under the flow, the process (6.11) is stationary. Since the points of $P$ are periodic for the flow, so are those of $P \setminus C_L$. Hence the stationary process (6.11) is periodic by Theorem 3.1. This shows that the process $X_\alpha^{C \setminus P}$ in the decomposition (6.6) has a nontrivial periodic component which contradicts Theorem 6.1. $\square$

COROLLARY 6.1. *Under the assumptions of Proposition* 6.1, *we also have $C_P = P$ a.e. $\nu(ds)$, where $P$ is the set of periodic points of the flow.*

COROLLARY 6.2. *A S$\alpha$S, $\alpha \in (0, 2)$, stationary process $X_\alpha$ with a minimal representation* (1.1) *is a stationary periodic (cyclic, resp.) process if and only if*

$$S = P \qquad \nu\text{-a.e.}\ (S = L\ \nu\text{-a.e., resp.}),$$

*where $P$ and $L$ are the periodic and the cyclic points of the generating flow, respectively. The equivalent is true if and only if the generating flow is periodic (cyclic, resp.).*

PROOF. Consider the case of stationary periodic processes. If $S = P$ $\nu$-a.e., then $S = C_P$ $\nu$-a.e. since, by Corollary 6.1, $C_P = P$ $\nu$-a.e. for minimal representations. Hence, by Theorem 4.1, the process $X_\alpha$ is a stationary periodic process. Conversely, if $X_\alpha$ is a stationary periodic process with a minimal representation (1.1), then $C_P = S$ $\nu$-a.e. by Theorem 4.1 and $C_P = P$ $\nu$-a.e. by Corollary 6.1. This implies $S = P$ $\nu$-a.e. The case of stationary cyclic processes can be considered in the same way by using Theorem 5.1 and Proposition 6.1. $\square$

Gathering the previous results, we obtain the following unique decomposition of S$\alpha$S stationary processes into four independent components.

THEOREM 6.2. *Let* $\{X_\alpha(t)\}_{t \in T}$ *be a S$\alpha$S, $\alpha \in (0,2)$, stationary process with a representation* (1.1). *Then the process* $X_\alpha$ *can be decomposed uniquely in distribution into four independent processes*

$$(6.12) \qquad X_\alpha \stackrel{d}{=} X_\alpha^{(1)} + X_\alpha^{(2)} + X_\alpha^{(3)} + X_\alpha^{(4)},$$

*where*

$$X_\alpha^{(1)}(t) = X_\alpha^D(t) = \int_D f_t(s) M_\alpha(ds),$$

$$X_\alpha^{(2)}(t) = X_\alpha^F(t) = \int_{C_F} f_t(s) M_\alpha(ds),$$

$$X_\alpha^{(3)}(t) = X_\alpha^L(t) = \int_{C_L} f_t(s) M_\alpha(ds),$$

$$X_\alpha^{(4)}(t) = X_\alpha^{C \setminus P}(t) = \int_{C \setminus C_P} f_t(s) M_\alpha(ds),$$

*and the sets $D$, $C$, $C_F$, $C_L$ and $C_P$ are defined in* (6.1), (6.2), (5.1), (5.3) *and* (4.1), *respectively. Here:*

1. *The process $X_\alpha^{(1)}$ has a mixed moving average representation* (1.7) *and is generated by a dissipative flow.*
2. *The process $X_\alpha^{(2)}$ is a harmonizable process with the representation* (1.9) *in the complex-valued case and is a trivial process with the representation* (1.11) *in the real-valued case.*
3. *The process $X_\alpha^{(3)}$ is a stationary cyclic process in the sense of Definition* 3.2.
4. *The process $X_\alpha^{(4)}$ is a stationary process generated by a conservative flow without a periodic component.*

*If the process $X_\alpha$ is generated by a flow $\{\phi_t\}_{t \in T}$, then the sets $D$ and $C$ are identical to the dissipative and the conservative parts of the flow $\{\phi_t\}_{t \in T}$,* *respectively. If, in addition, the representation of the process $X_\alpha$ is minimal, then the sets $C_F$, $C_L$ and $C_P$ are the fixed, cyclic and periodic points of the flow $\{\phi_t\}_{t \in T}$, respectively.*

PROOF. The theorem follows from the decomposition of a S$\alpha$S stationary process into three components in Rosiński (1995), Theorem 6.1, Proposition 6.1 and Corollary 6.1. □

## APPENDIX A

**Flows on a standard Lebesgue space.** We provide here a number of definitions related to flows which are used throughout the paper. A measure space $(S, \mathcal{S}, \nu)$ is called a *standard Lebesgue space* when $(S, \mathcal{S})$ is a standard Borel space equipped with a $\sigma$-finite measure $\nu$. A standard Borel space is a measurable space measurably isomorphic (i.e., there is a one-to-one, onto and bimeasurable map) to a Borel subset of a complete separable metric space. We may thus suppose without loss of generality that a standard Borel space is a subset of a complete separable metric space. The corresponding $\sigma$-field $\mathcal{S}$ is defined as the smallest $\sigma$-field that contains all Borel sets. Standard Lebesgue spaces (or standard Borel spaces) are convenient to work with, have nice properties and are widely used in ergodic theory [see Walters (1982) and Petersen (1983)] and in other areas of mathematics [see Zimmer (1984), Arveson (1976) and Mackey (1957)]. They were used by Rosiński (1995) and Pipiras and Taqqu (2002a, b) in the context of stable processes. The Euclidean space equipped with a Lebesgue measure, for example, is a standard Lebesgue space.

A flow $\{\phi_t\}_{t \in T}$ with $T = \mathbb{R}$ or $T = \mathbb{Z}$ on a standard Lebesgue space $(S, \mathcal{S}, \nu)$ is a collection of deterministic maps $\phi_t : (S, \mathcal{S}) \mapsto (S, \mathcal{S})$ such that $\phi_0(s) = s$ and

(A.1) $\qquad \phi_{t_1+t_2}(s) = \phi_{t_1}(\phi_{t_2}(s)) \qquad$ for all $t_1, t_2 \in T, s \in S.$

A flow $\{\phi_t\}_{t \in T}$ is called *nonsingular* if, for every $t \in T$, $\nu(N) = 0$ if and only if $\nu(\phi_t^{-1}(N)) = 0$. It is called *measurable* if the map $\phi_t(s) : T \times S \mapsto S$ is measurable.

It is known that a measurable nonsingular flow on a space $(S, \mathcal{S}, \nu)$ has the so-called *Hopf decomposition* [see Krengel (1985), page 17, Rosiński (1995), page 1171, or Pipiras and Taqqu (2002a, b)]. The Hopf decomposition is a (a.e.) partition of the space $S$ into two disjoint sets $C$ and $D$. The set $D$, called a *dissipative* part of the flow, can be written as $D = \sum_{k \in \mathbb{Z}} \phi_1^k(B)$ for some wandering set $B$. ["Wandering" means that the sets $\phi_1^m(B)$ and $\phi_1^n(B)$ are disjoint for $m \neq n$.] The set $C$, called a *conservative* part of the flow, is such that it has no wandering set of positive measure. Moreover, the sets $C$ and $D$ can be taken to be invariant under the flow [i.e., $\phi_t^{-1}(C) =$ $C$ and $\phi_t^{-1}(D) = D$ for all $t \in T$]. The flow $\{\phi_t\}_{t \in T}$ is called conservative (dissipative, resp.) if $S = C$ ($S = D$, resp.) a.e. For a general flow $\{\phi_t\}_{t \in T}$, its restriction $\{\phi_t|_C\}_{t \in T}$ ($\{\phi_t|_D\}_{t \in T}$, resp.) is a conservative (dissipative, resp.) flow.

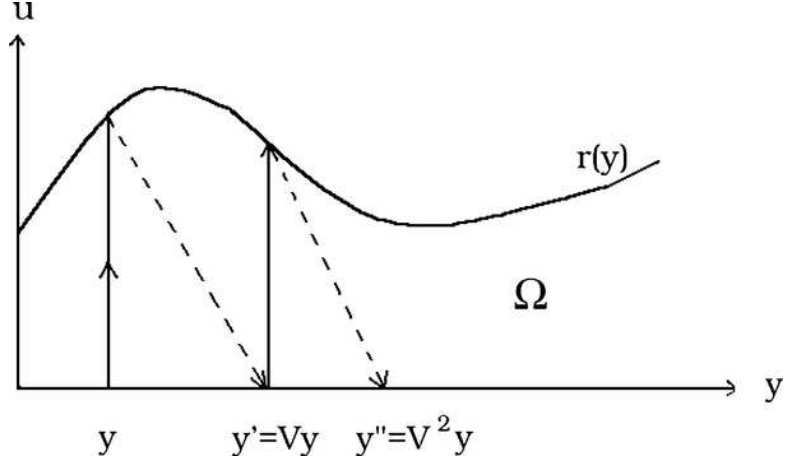

Fig. 1. *View the flow $\widetilde{\phi}_t$ as moving up vertically at constant speed until it reaches the level $r(y)$ and then jumps back to a point $(y',0)$ before it renews its vertical climb, this time from the point $y'$.*

In this work, we use the notion of a *special flow* $\{\widetilde{\phi}_t\}_{t \in T}$. Informally, the flow $\widetilde{\phi}_t(y, u)$ is defined on the set of points

$$\Omega = \{(y, u) : 0 \leq u < r(y), y \in Y\} = Y \times [0, r(\cdot)),$$

where $r(y)$ is a positive function. Plotting $(y, u)$ in two dimensions, we can view the flow $\widetilde{\phi}_t$ as moving up vertically at constant speed until it reaches the level $r(y)$, and then jumps back to a point $(y', 0)$ before it renews its vertical climb, this time from the point $y'$ (see Figure 1). Thus, if we focus only on the horizontal $Y$ axis, the flow starting at $y$ moves to $y' = Vy$, then to $V^2 y, \ldots, V^n y, \ldots$. Since the flow $\widetilde{\phi}_t$ moves constantly, observe that it has no fixed points.

The flow $\widetilde{\phi}_t$ is defined formally as follows. Let $(Y, \mathcal{Y}, \tau)$ be a standard Lebesgue space, let $V$ be a one-to-one, onto, bimeasurable and nonsingular map of $Y$ onto itself, and let $r$ be a positive measurable function on $Y$ such that $\sum_{k=0}^{\infty} r(V^k y) = \sum_{k=-\infty}^{-1} r(V^k y) = \infty$. Set $\Omega = \{(y, u) : 0 \leq u < r(y), y \in Y\}$, $\mathcal{E} = \mathcal{Y} \otimes \mathcal{B}([0, r(\cdot)))$ and let $P$ be a measure on $(\Omega, \mathcal{E})$ such that $dP(y, u) = p(y, u)\tau(dy)\,du$ and $P(\Omega) = 1$. Consider now the map defined on $\Omega$ by

$$(A.2) \qquad \widetilde{\phi}_t(y, u) = (V^n y, u + t - r_n(y))$$

for

$$0 \le u + t - r_n(y) < r(V^n y),$$

where $r_n(y) = \sum_{k=0}^{n-1} r(V^k y)$ if $n \ge 1$, $r_n(y) = 0$ if $n = 0$, and $r_n(y) = \sum_{k=n}^{-1} r(V^k y)$ if $n \le -1$. We can verify that $\{\widetilde{\phi}_t\}_{t \in \mathbb{R}}$ is a (measurable, nonsingular) flow on $(\Omega, \mathcal{E}, P)$. It is called a special flow built under the function $r$. According to Theorem 3.1 in Kubo (1969), a (measurable, nonsingular) flow $\{\phi_t\}_{t \in \mathbb{R}}$ without fixed points on a standard Lebesgue space is null isomorphic (mod 0) to some special flow $\{\widetilde{\phi}_t\}_{t \in \mathbb{R}}$ built under the function $r$.

In addition to flows, we also use a related functional called a *cocycle*. Let $A$ be a second countable group, that is, a topological group that has a countable base for the topology. For example, $A = \{-1, 1\}$ or $A = \{w : |w| = 1\}$ with a multiplication operation and $A = \mathbb{R}$ with an addition operation. A measurable map $a_t(s) : T \times S \mapsto A$ is called a cocycle for a measurable flow $\{\phi_t\}_{t \in T}$ if

(A.3) $$a_{t_1+t_2}(s) = a_{t_1}(s) a_{t_2}(\phi_{t_1}(s)) \qquad \forall t_1, t_2 \in T, s \in S.$$

In this paper, we use exclusively the cases $A = \{-1, 1\}$ and $A = \{w : |w| = 1\}$, but cocycles are typically associated in the literature [see Zimmer (1984)] with second countable groups.

## APPENDIX B

**Minimal representations for stable processes.** Finally, we define minimal integral representations of stable processes which play a central role in relating stable processes to flows. An integral representation $\{f_t\}_{t \in T} \subset L^\alpha(S, \mathcal{S}, \nu)$ is called minimal for the process $X_\alpha$ given by (1.1) if:

(M1) $\mathrm{supp}\{f_t, t \in T\} = S \quad \nu$-a.e.;
(M2) $\sigma\{f_u/f_v, u, v \in T\} = \mathcal{S}$ modulo $\nu$

[see Hardin (1982), Rosiński (1995, 1998a) and Pipiras and Taqqu (2002a)]. A condition equivalent to (M2) is the following:

(M2′) For every nonsingular map $\phi : S \mapsto S$ and $h : S \mapsto \mathbb{R} \setminus \{0\}$ such that, for each $t \in T$,

(B.1) $$f_t(s) = h(s) f_t(\phi(s)) \qquad \text{a.e. } \nu(ds),$$

we have $\phi(s) = s$ a.e. $\nu(ds)$ [see Rosiński (1998a)].

As shown in Hardin (1982), every separable in probability S$\alpha$S process has a minimal integral representation. Rosiński (1995) showed that minimal integral representations of stationary S$\alpha$S processes are related to flows in the sense of Definition 1.1.

**Acknowledgments.** We thank an anonymous referee for carefully reading the original manuscript and for numerous suggestions.

## REFERENCES

ARVESON, W. (1976). *An Invitation to $C^*$-Algebras*. Springer, New York.

BINGHAM, N. H., GOLDIE, C. M. and TEUGELS, J. L. (1987). *Regular Variation*. Cambridge Univ. Press. MR512360

COHN, D. L. (1972). Measurable choice of limit points and the existence of separable and measurable processes. *Z. Wahrsch. Verw. Gebiete* **22** 161–165. MR305444

GROSS, A. (1994). Some mixing conditions for stationary symmetric stable stochastic processes. *Stochastic Process. Appl.* **51** 277–295. MR1288293

HALMOS, P. R. (1950). *Measure Theory*. Van Nostrand, New York. MR33869

HARDIN, C. D., JR. (1982). On the spectral representation of symmetric stable processes. *J. Multivariate Anal.* **12** 385–401. MR666013

KRENGEL, U. (1985). *Ergodic Theorems*. de Gruyter, Berlin. MR797411

KUBO, I. (1969). Quasi-flows. *Nagoya Math. J.* **35** 1–30. MR247032

MACKEY, G. W. (1957). Borel structure in groups and their duals. *Trans. Amer. Math. Soc.* **85** 134–165. MR89999

PETERSEN, K. (1983). *Ergodic Theory*. Cambridge Univ. Press.

PIPIRAS, V. and TAQQU, M. S. (2002a). Decomposition of self-similar stable mixed moving averages. *Probab. Theory Related Fields* **123** 412–452. MR833286

PIPIRAS, V. and TAQQU, M. S. (2002b). The structure of self-similar stable mixed moving averages. *Ann. Probab.* **30** 898–932. MR1905860

ROSIŃSKI, J. (1995). On the structure of stationary stable processes. *Ann. Probab.* **23** 1163–1187. MR1349166

ROSIŃSKI, J. (1998a). Minimal integral representations of stable processes. Preprint.

ROSIŃSKI, J. (1998b). Structure of stationary stable processes. In *A Practical Guide to Heavy Tails: Statistical Techniques and Applications* (R. Adler, R. Feldman and M. S. Taqqu, eds.) 461–472. Birkhäuser, Boston. MR1652294

ROSIŃSKI, J. (2000). Decomposition of stationary $\alpha$-stable random fields. *Ann. Probab.* **28** 1797–1813.

ROSIŃSKI, J. and SAMORODNITSKY, G. (1996). Classes of mixing stable processes. *Bernoulli* **2** 365–377. MR1813849

ROZANOV, Y. A. (1967). *Stationary Random Processes*. Holden-Day, San Francisco.

SAMORODNITSKY, G. and TAQQU, M. S. (1994). *Stable Non-Gaussian Processes: Stochastic Models with Infinite Variance*. Chapman and Hall, New York.

SURGAILIS, D., ROSIŃSKI, J., MANDREKAR, V. and CAMBANIS, S. (1998). On the mixing structure of stationary increment and self-similar S$\alpha$S processes. Preprint. MR214134

WALTERS, P. (1982). *An Introduction to Ergodic Theory*. Springer, New York. MR648108

ZIMMER, R. J. (1984). *Ergodic Theory and Semisimple Groups*. Birkhäuser, Boston. MR776417

DEPARTMENT OF MATHEMATICS
BOSTON UNIVERSITY
111 CUMMINGTON STREET
BOSTON, MASSACHUSETTS 02215
USA
E-MAIL: pipiras@math.bu.edu
E-MAIL: murad@math.bu.edu